\documentclass{paper}

\usepackage{arxiv}

\usepackage[utf8]{inputenc} % allow utf-8 input
\usepackage[T1]{fontenc}    % use 8-bit T1 fonts
\usepackage{hyperref}       % hyperlinks
\usepackage{float}			% floating pics and tabulars
\usepackage{url}            % simple URL typesetting
\usepackage{booktabs}       % professional-quality tables
\usepackage{amsfonts}       % blackboard math symbols
\usepackage{nicefrac}       % compact symbols for 1/2, etc.
\usepackage{microtype}      % microtypography
\usepackage{mathrsfs,bm,siunitx}
\usepackage{mathtools}
\usepackage{ragged2e}
\usepackage[]{graphicx}
\usepackage[]{color}
\usepackage{alltt}
\usepackage{mathtools}
\usepackage{amssymb}
\usepackage{aligned-overset}
\usepackage{amsmath}
\numberwithin{equation}{section}
\usepackage{tikz}
\usepackage{capt-of}
\usepackage{xcolor,pgf}
\usepackage{color}
\usepackage{transparent,mathtools}
\usepackage{caption}
\usetikzlibrary{arrows}
\usetikzlibrary{decorations.markings}
\usepackage{pgfplots}
    \usetikzlibrary{
        angles,
        quotes,
    }
\usepackage[square, comma, numbers, sort]{natbib}
\usepackage{subfig}

\usepackage{lipsum}
\usepackage{graphicx}
\usepackage{amsmath,amsthm}
\numberwithin{equation}{section}
\usepackage{color}
\usepackage{colortbl}
\usepackage{xcolor}
\colorlet{LighterGray}{white!70!gray}
\colorlet{LightererGray}{white!85!gray}
\usepackage{transparent}
\graphicspath{{Bilder/}}
\usepackage{tikz}
\usetikzlibrary{arrows}
\usetikzlibrary{decorations.markings}
\usepackage{pgfplots}
\pgfplotsset{compat=1.14}
    \usetikzlibrary{
        angles,
        quotes,
    }
    
    \usepackage{multirow}
\usepackage[utf8]{inputenc}
\usepackage{ragged2e}
\usepackage{xcolor}
\colorlet{LighterGray}{white!70!gray}
\colorlet{LightererGray}{white!85!gray}
\usepackage[square, comma, numbers, sort]{natbib}
\usepackage{fancyhdr}
\newcommand{\changefont}{%
    \fontsize{8}{11}\selectfont
}
\usepackage{hyperref}
\hypersetup{
  colorlinks = true,
  linkcolor = magenta,
  urlcolor = magenta,
  citecolor = cyan}
\usepackage{scalerel,stackengine,amsmath}

\usepackage{enumitem}
\newlist{abbrv}{itemize}{1}
\setlist[abbrv,1]{label=,labelwidth=1in,align=parleft,itemsep=0.1\baselineskip,leftmargin=!}
\usepackage{listings}
\usepackage{color} %red, green, blue, yellow, cyan, magenta, black, white
\definecolor{mygreen}{RGB}{28,172,0} % color values Red, Green, Blue
\definecolor{mylilas}{RGB}{170,55,241}

\usepackage{amstext} % for \text macro
\usepackage{array}   % for \newcolumntype macro
\newcolumntype{L}{>{$}l<{$}}
\usepackage{bm}
\DeclareMathOperator{\logg}{logg}
\DeclareMathOperator{\l2}{\ell^2(\mathbb{N})}

\usepackage{pgfplots}
\pgfplotsset{compat=1.14}
\newtheorem{theorem}{Theorem}[section]

\theoremstyle{definition}
\newtheorem{dfn}{Definition}[section]

\newtheorem{assump}{Assumption}[section]
\newtheorem{alg}{Algorithm}[section]

\title{Bayesian identification of material parameters in viscoelastic structures as an inverse problem in a semigroup setting}

\author{
  Rebecca Rothermel \\
  Department of Mathematics\\
  Saarland University\\
  Saarbr\"ucken, Germany \\
  \texttt{klein@num.uni-sb.de} \\
       \And
  Thomas Schuster \\
  Department of Mathematics\\
  Saarland University\\
  Saarbr\"ucken, Germany \\
  \texttt{thomas.schuster@num.uni-sb.de} \\
}

\begin{document}
\maketitle

\begin{abstract}

The article considers the nonlinear inverse problem of identifying the material parameters in viscoelastic structures based on a generalized Maxwell model. The aim is to reconstruct the model parameters from stress data acquired from a relaxation experiment, where the number of Maxwell elements, and thus the number of material parameters themselves, are assumed to be unknown. This implies that the forward operator acts on a Cartesian product of a semigroup (of integers) and a Hilbert space and demands for an extension of existing regularization theory. We develop a stable reconstruction procedure by applying Bayesian inversion to this setting. We use an appropriate binomial prior which takes the integer setting for the number of Maxwell elements into account and at the same time computes the underlying material parameters. We extend the regularization theory for inverse problems to this special setup and prove existence, stability and convergence of the computed solution. The theoretical results are evaluated by extensive numerical tests.
%Thus, the number of material parameters depends on the solution of the inverse problem. 

\end{abstract}

% keywords can be removed
\keywords{viscoelastic material, Bayesian inversion, semigroup, Maxwell model, binomial prior, inverse problem}  %TODO
%\newpage

%%%%%%%%%%%%%%%%%%%%%%%%

\section{Introduction}\label{sec:intro}
%(mit kurzer Problembeschreibung sowie Literaturüberblick)

%%%%%%%%%%%%%%%%%%%%%%%%
We consider the problem of identifying material parameters in viscoelastic structures. Parameter identification represents a challenging class of inverse problems having important and demanding real-world applications. These include, for example, the identification of the distortion energy density of hyperelastic materials \cite{Hartmann2003, Woestehoff2015, Binder2015, Seydel2017OnField, Seydel2017IdentifyingMethod, Klein2021}, the surface enthalpy-dependent heat fluxes of steel plates \cite{Rothermel2021, rothermel2019parameter}, inverse scattering problems \cite{colton1998inverse}, the estimation of parameters from waveform information \cite{eller2024tangential}, the inverse kinematic problem \cite{klibanov2016reconstruction}, electrical impedance tomography \cite{borcea2002electrical} or terahertz tomography \cite{Wald2018Tomographic} to name only a few. All these problems have in common that they are usually non-linear and (locally) ill-posed. That means, that even small errors in the measured data lead to large inaccuracies in the computed solution, if no regularization is applied. There exists a vast amount of literature for solving inverse problems in Hilbert and Banach spaces such as, e.g., \cite{Engl1996RegularizationProblems, Louis1989InverseProbleme, Rieder2003KeineProblemen, Kirsch2011, Schuster2012RegularizationSpaces, Kaltenbacher2008IterativeProblems}. The problem which is considered in the presented article needs for an extension of the existing theory to spaces which show only the structure of a semigroup.

Understanding how a material deforms when a force is applied is essential for many industrial applications ranging from food processing, additive manufacturing, structural health monitoring \cite{giurgiutiu2007structural}, seismics \cite{rieder2020all} to product design \cite{Bonfanti2020FractionalMaterials}.

%Eventuell nur die Referenzen hier raus und nach oben rein
In food processing, for example, the texture of bread depends very much on its mechanical properties, which in turn are determined by the manufacturing process \cite{lazaridou2007effects, tanner2008bread}.
In additive manufacturing, plastics are melted in the form of polymer filaments and then built up layer by layer to form a three-dimensional workpiece. The flow behavior of the material uses both affects, the processing time and the strength of the printed object \cite{mackay2018importance, corker20193d}. The mechanical response of polymers can be used to determine the formation of their microstructure \cite{hata1968effect, van2010decoding}.

The scientific study of material deformation is called \emph{rheology}. There are three basic models that describe the idealized behavior of materials: elasticity, viscosity and plasticity. In reality, however, most materials show a combination of different properties that occur at various degrees \cite{meinhard1999rheologische}.
In this work, we consider viscoelastic materials that exhibit both, viscous and elastic properties under an applied force. Polymers represent a typical example for such materials. A viscoelastic model should be capable to describe two processes: \emph{relaxation} and \emph{creep}. Relaxation describes the process, during which viscoelastic materials relax while deformation remains constant, i.e., the stress in the material decreases. If a force is applied to a viscoelastic material, the deformation occurs with a time delay what is referred to as creep. This is discussed in section \ref{chap:Visko}. Further theory on the rheology and material behavior of viscoelastic materials can be found in \cite{Tschoegl1989TheBehavior, Wineman2000MechanicalPolymers}. Unlike for elastic materials, for viscoelastic materials the constitutive equations are often unknown and there exists a variety of different models to describe viscoelasticity.

There are several publications dealing with the identification of viscoelastic  parameters. However, these differ from the present work by several features. 
In \cite{Shukla2021,Emri1998} viscoelastic structures are considered and their material parameters identified, but instead of the generalized Maxwell model they use a different modeling. A comparable model is used in \cite{Diebels2018}. There, however, only the stiffnesses of a viscoelastic material are reconstructed. The relaxation times as well as the number of Maxwell elements are assumed to be known. 
Relaxation experiments in combination with cyclic tests are used in \cite{Scheffer2013Optimisationmaterials} to determine the basic stiffness and the material parameters of four Maxwell elements. Here, the number of Maxwell elements is known a priori, see also \cite{Fernanda2011,Park1999, Sharma2020MoistureProperties}. Babaei et al. \cite{babaei2016efficient} propose two methods for solving the underlying inverse problem. The so-called \emph{ad hoc method} first guesses the number of Maxwell elements $n$ and subsequently reconstructs the material parameters for fixed $n$. Our approach differs from this procedure since the computation of $n$ is part of the algorithm and the computation of the parameters is done simultaneously.
The second approach (\emph{discrete spectral approach}) in \cite{babaei2016efficient} distributes about 1000 relaxation times equidistantly in a logarithmic scale over the interval $[10^{-1},10^3]$. This drastically simplifies the reconstruction and only requires solving a system of linear equations to calculate the stiffnesses. Then, the dominant parameters are identified and the correct number of Maxwell elements is found, although this process and the handling of the remaining parameters are not described in more detail. 

We use a generalized Maxwell model, which is characterized by an unknown number of Maxwell elements $n$ and material parameters (relaxation times $\tau_j$ and stiffnesses $(\mu, \mu_j$). This means, that the exact solution determines at the same time the number of parameters to be determined, a feature that has to be included in the modeling process. The forward operator, whose construction is outlined in section \ref{sec: Vorwärtsoperator}, maps the material parameters and the number of Maxwell elements to the stress function. This function describes the time history of stress in the material during a relaxation experiment, where a strain is applied to the material and kept constant. The forward operator acts on a Cartesian product of a semigroup (the integers $\mathbb{N}$) and a Hilbert space ($\l2$).
In section \ref{sec: inverse problem} we define the inverse problem of determining the number of Maxwell elements and simultaneously the material parameters from the stress function.\\

Unfortunately, a large number of parameters often leads to overfitting a noisy data term and thus to unavoidable errors in the parameters. 
In a previous article the authors in  \cite{Rothermel2022} developed a clustering algorithm adapted to this problem. But it shows a strong error susceptibility to noisy data.
As a solution, we propose in section \ref{chap: Bayes Inversion} a novel method using statistical Bayesian inversion theory. This uses a binomial prior to estimate the number of Maxwell elements. The deduction of the method is subject of section \ref{sec: Bayes for us}.
The developed algorithm alternately searches for a suitable solution for $n$ and the material parameters by minimizing an appropriate Tikhonov functional. the minimization is done in $\mathbb{N}\times \l2$ and, hence, the standard theory for regularizing Tikhonov functionals does not apply, since it uses Hilbert and Banach spaces settings (c.f. \cite{Engl1996RegularizationProblems, Schuster2012RegularizationSpaces, Ito2015InverseAlgorithms, Hofmann2007AOperators}). To this end we extend the reguarization theory for Tikhonov functionals to this particular setting where we consider the integers $\mathbb{N}$ as topological semigroup endowed with the discrete topology. Proofs for convergence and stability of our proposed algorithm is subject of section \ref{sec:regularization_algorithm}. 

Numerical validation of the theoretical framework is done in section \ref{chap:Numeric}. For this purpose, we introduce the clustering algorithm from \cite{Rothermel2022} in order to compare the reconstruction results of the different algorithms. We perform experiments using different exact and noise perurbed data sets. We analyze different displacement rates of the strain function in the relaxation experiment and consider the effect of the success probability associated with the binomial distribution of the prior. Additionally, we introduce different penalty terms with respect to the material parameters and analyze their influence on the reconstruction results.\\

Summarizing, the article contains the following innovations:

\begin{itemize}
    \item extending regularization theory to forward operators acting on the Cartesian product $\mathbb{N} \times \l2$, where the integers form a topological semigroup,
    \item dependence of the number of material parameters on the number of Maxwell elements $n$, and thus on part of the inverse problem solution itself,
    \item Bayesian inversion approach using a binomial prior,
    \item proof of regularization property (existence, stability, convergence of solutions) for our approach
\end{itemize}

%%%%%%%%%%%%%%%%%%%%%%%
\section{A rheological model for viscoelastic material behavior} \label{chap:Visko}
%mit Beschreibung des rheologischen Modells und des daraus abgeleiteten Vorwärtsproblems

%%%%%%%%%%%%%%%%%%%%%%%

%%%%%%%%%%%%%%%%%%%%%%%%%%%%%%%%%%%%%%%%%%%%%%%%%%

\subsection{Viscoelastic materials}
\label{sec_visco}

We introduce the rheological model of a viscoelastic material, which we use as starting point. 
%We will determine the material parameters of a viscoelastic material by measuring the stress in a relaxation experiment. 
For comprehensive introductions to the phenomenological behavior and modeling of viscoelasticity we refer to textbooks such as \cite{Tschoegl1989TheBehavior, Wineman2000MechanicalPolymers}.

A viscoelastic material is characterized by a combination of viscous as well as elastic behavior. Viscosity means that, if a force is applied, deformations are observed that are unlimited and irreversible. In contrast, elasticity describes the ability of a body to deform back to its original shape after the applied force is removed. In this case the deformations are limited reversible. Viscoelastic behavior is exhibited by various polymeric materials such as adhesives, elastomers and rubber. 

\begin{figure}[h]
\centering
\def\svgwidth{.6\linewidth}
%% Creator: Inkscape 1.0 (4035a4f, 2020-05-01), www.inkscape.org
%% Accompanies image file 'StrainsExplained.pdf' (pdf, eps, ps)
\begingroup%
  \makeatletter%
  \providecommand\color[2][]{%
    \errmessage{(Inkscape) Color is used for the text in Inkscape, but the package 'color.sty' is not loaded}%
    \renewcommand\color[2][]{}%
  }%
  \providecommand\transparent[1]{%
    \errmessage{(Inkscape) Transparency is used (non-zero) for the text in Inkscape, but the package 'transparent.sty' is not loaded}%
    \renewcommand\transparent[1]{}%
  }%
  \providecommand\rotatebox[2]{#2}%
  \newcommand*\fsize{\dimexpr\f@size pt\relax}%
  \newcommand*\lineheight[1]{\fontsize{\fsize}{#1\fsize}\selectfont}%
  \ifx\svgwidth\undefined%
    \setlength{\unitlength}{650.25bp}%
    \ifx\svgscale\undefined%
      \relax%
    \else%
      \setlength{\unitlength}{\unitlength * \real{\svgscale}}%
    \fi%
  \else%
    \setlength{\unitlength}{\svgwidth}%
  \fi%
  \global\let\svgwidth\undefined%
  \global\let\svgscale\undefined%
  \makeatother%
  \begin{picture}(1,0.69665513)%
    \lineheight{1}%
    \setlength\tabcolsep{0pt}%
    \put(0,0){\includegraphics[width=\unitlength,page=1]{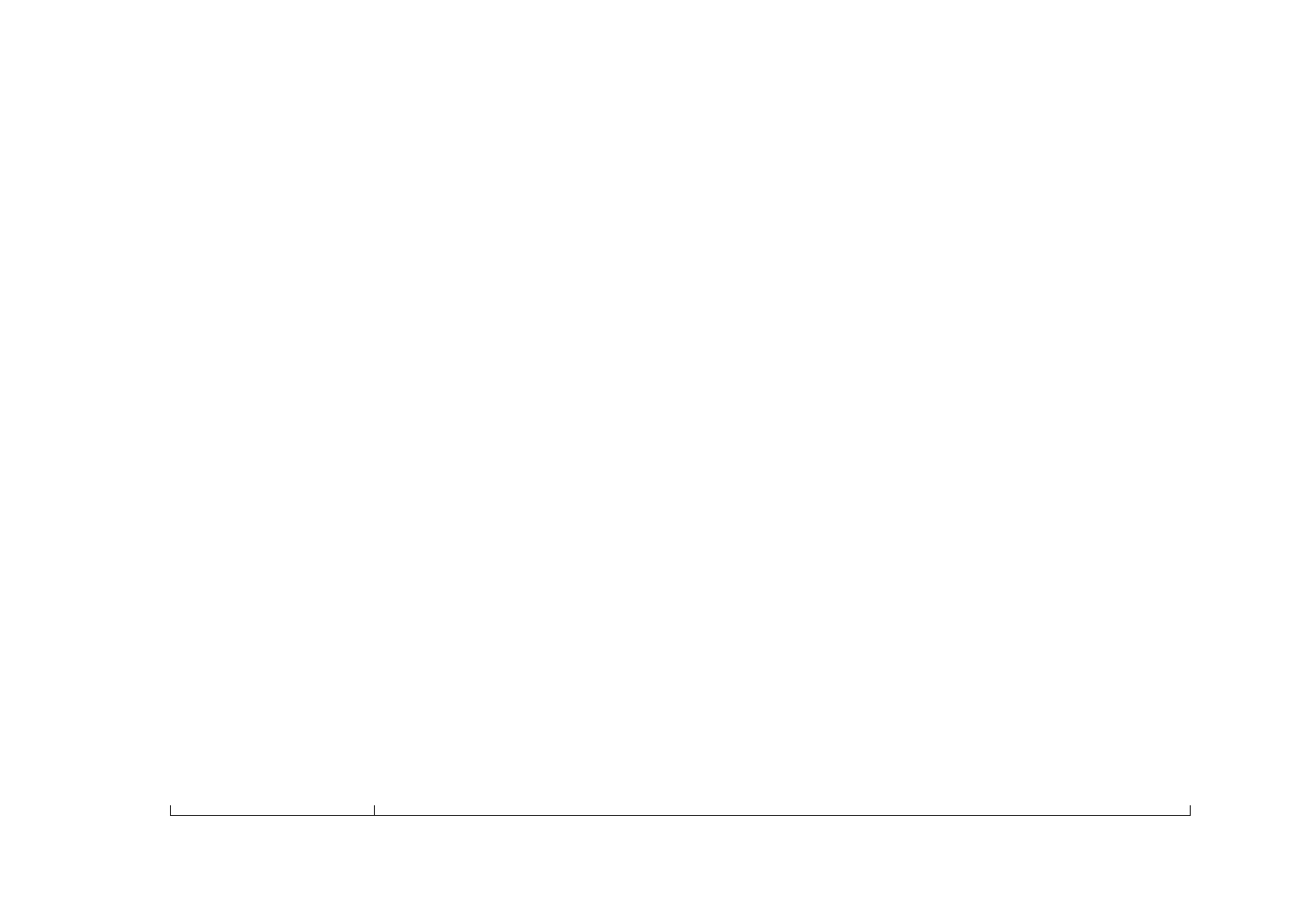}}%
    \put(0.49877357,0.02150159){\makebox(0,0)[lt]{\lineheight{1.25}\smash{\begin{tabular}[t]{l}Time [s]\end{tabular}}}}%
    \put(0,0){\includegraphics[width=\unitlength,page=2]{StrainsExplained.pdf}}%
    \put(0.09108519,0.25445268){\rotatebox{90}{\makebox(0,0)[lt]{\lineheight{1.25}\smash{\begin{tabular}[t]{l}Strain \end{tabular}}}}}%
    \put(0.09568491,0.61126243){\makebox(0,0)[lt]{\lineheight{1.25}\smash{\begin{tabular}[t]{l}$\bar{\varepsilon}$\end{tabular}}}}%
    \put(0.10688273,0.03487846){\makebox(0,0)[lt]{\lineheight{1.25}\smash{\begin{tabular}[t]{l}$t=0$\end{tabular}}}}%
    \put(0.24507454,0.03487846){\makebox(0,0)[lt]{\lineheight{1.25}\smash{\begin{tabular}[t]{l}$t=\frac{\bar{\varepsilon}}{\dot{\varepsilon}}$\end{tabular}}}}%
    \put(0.8906896,0.03487784){\makebox(0,0)[lt]{\lineheight{1.25}\smash{\begin{tabular}[t]{l}$t=T$\end{tabular}}}}%
  \end{picture}%
\endgroup%

\caption{Strain curve $\varepsilon(t)$ with strain rate $\dot{\varepsilon}$ and maximum strain value $\bar{\varepsilon}$ }
\label{abb:strain_exp}
\end{figure}

\begin{figure}[h]
\centering
\def\svgwidth{.6\linewidth}
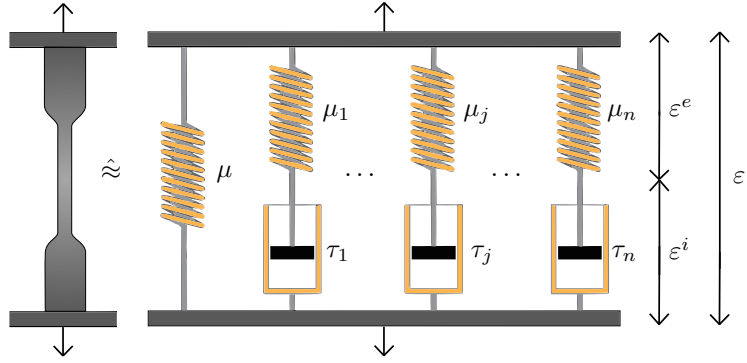
\caption{Standard specimen and generalized Maxwell model with unknown number of Maxwell elements $n$ and $2n+1$ material parameters}
\label{abb:rheologicalModel}
\end{figure}

Figure \ref{abb:rheologicalModel} shows the typical shape of a standard specimen which is clamped at the thick ends and loaded in the direction of the arrows. According to Saint Venant's principle, the disturbances caused by the clampings at the ends of the specimen decay after a short distance. Therefore, we can assume that the strain and stress state at the center of the specimen is homogeneous and we can use a one-dimensional model. The force and extension of the specimen are measured and can be used to make direct calculations from strain and stress. 
We consider the following relaxation experiment: 
For a given strain rate $\dot{\varepsilon}$ and a maximum strain value $\bar{\varepsilon}$, the strain in a time interval $t \in [0,T]$ is given as
\begin{align}
\varepsilon(t)= \left\{\begin{array}{ll} \dot{\varepsilon} \cdot t, &0 \leq t\leq \frac{\bar{\varepsilon}}{\dot{\varepsilon}} \\[1ex]
         \bar{\varepsilon}, & \frac{\bar{\varepsilon}}{\dot{\varepsilon}} <t \leq T. \end{array}\right. 
        \label{eqn:strain}
\end{align}
The function is plotted in figure \ref{abb:strain_exp} and describes the following procedure. The material is stretched until a maximum strain value $\bar{\varepsilon}$ is reached at a strain rate $\dot{\varepsilon}$. The strain rate is $\dot{\varepsilon}=\frac{\dot{\varepsilon}_u}{l_0}$ with the displacement rate $ \dot{\varepsilon}_u $, given by the testing device, and specimen length $l_0$.  Thus, the maximum strain is reached at time $t=\frac{\bar{\varepsilon}}{\dot{\varepsilon}}$. After that, the applied strain is kept constant. 
\\

\begin{figure}[h]
     \centering
     \subfloat{\def\svgwidth{0.5\textwidth}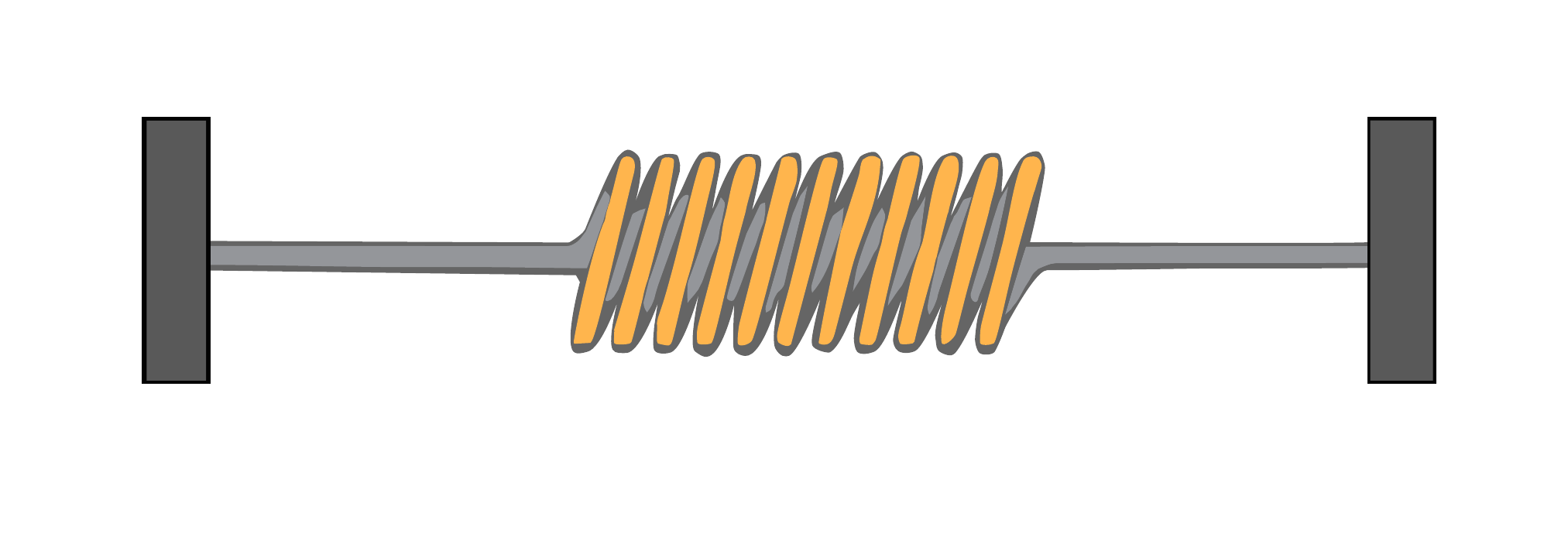} \\
     \subfloat{\def\svgwidth{0.5\textwidth}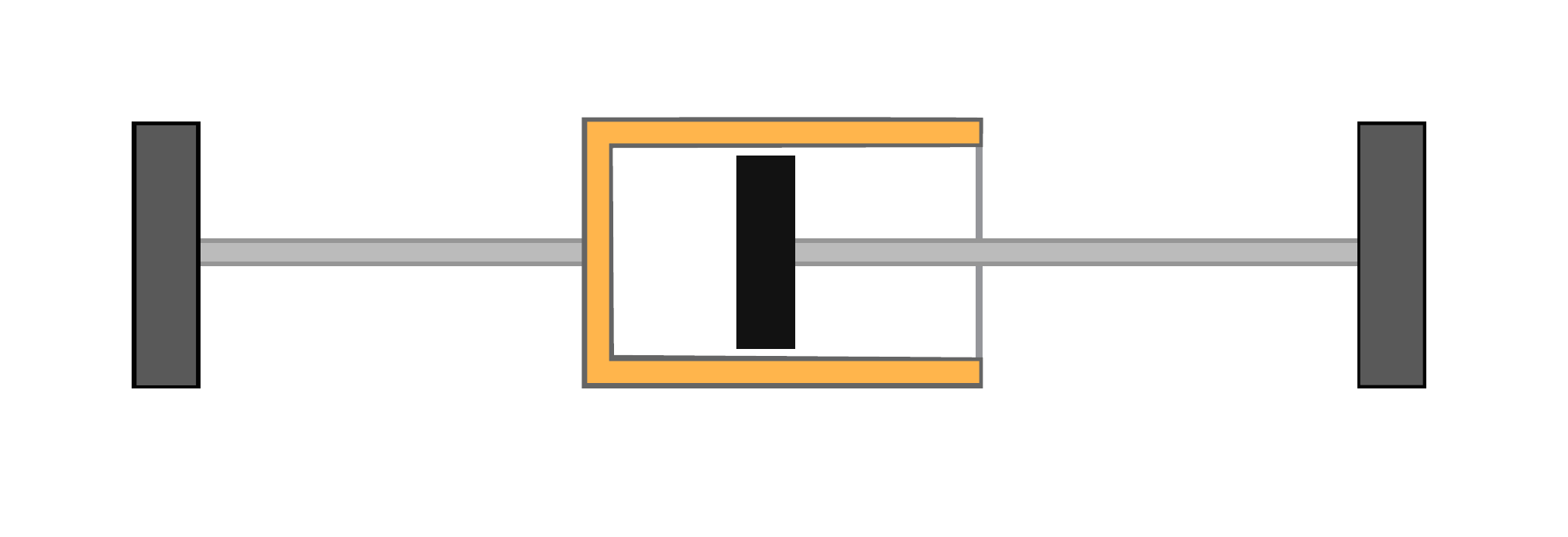}
\caption{A Hooke and a Newtonian element}
\label{abb:FederDämpfer}
\end{figure}

The stress in a viscoelastic material is characterized by a multi-parameter rheological model. 
In rheology, linear elastic behavior is modeled by a spring, also called a \emph{Hooke element} (see figure \ref{abb:FederDämpfer}). Hooke's law assumes that stress and strain depend linearly with the modulus of elasticity as proportional constant. Alternatively, Hooke's law can be used to link the displacement and the force using the spring's stiffness. The spring reacts immediately to the applied strain, but the deformation is limited and, if the force is released, the spring returns to its initial position like an ideal elastic body.
The second element, that is relevant for our model, is a damper or \emph{Newtonian element}, also shown in figure \ref{abb:FederDämpfer}. While for the Hooke body the measured force is proportional to the displacement, for the Newton body it is proportional to the velocity. The constant of proportionality is called the damping constant and corresponds to the viscosity.
This model for viscous material behavior responds to force with a time delay. In contrast to an elastic body, the deformation is unlimited as long as the force is applied. In this case, the deformation is irreversible and the deformation rate is proportional to the stress.
The parallel connection of a spring and a damper is called a \emph{Kelvin-Voigt element}. If a strain is applied to such a body and kept constant, the spring deforms elastically but is slowed down by the damper so that the deformation occurs with a time delay. This is also known as creep and is typical of a viscoelastic material.
However, this composition can only describe creep processes, but not relaxation processes. The latter can be implemented by a Maxwell element, the serial composition of a spring and a damper. If a force acts on a Maxwell body, it reacts with a rate of deformation that depends on the magnitude of the force as well as its rate of change. However, if the applied strain is kept constant, as in the relaxation experiment, the stress in the Maxwell body relaxes to zero, since the damper relaxes according to its relaxation time, which is given by the ratio of viscosity and stiffness. Thus, a simple Maxwell body represents a viscoelastic fluid, but cannot model a viscoelastic solid.\\

The simplest model to model a viscoelastic solid is a three-parameter model consisting of a parallel combination of a Maxwell element with a spring (compare \cite{Goldschmidt2015ModellingMoisture, Johlitz2007ExperimentalSystems, Scheffer2013Optimisationmaterials}).
If one applies strain to such a solid model, both springs will stretch. If the strain is then kept constant during the rest of the experiment, the damper expands according to its relaxation time.  When it is fully relaxed, the spring in the Maxwell element reaches its original undeformed position and the stress is fully absorbed by the single spring. 
this is why the model with a single Maxwell element can only simulate a relaxation time for the entire duration of the relaxation of the material. However, viscoelastic materials exhibit relaxation behavior with different relaxation times \cite{Sharma2020MoistureProperties, Goldschmidt2015ModellingMoisture}. Hence, we need to cover the entire relaxation spectrum (compare \cite{Scheffer2013Optimisationmaterials, Bonfanti2020FractionalMaterials, Jalocha2020PayneTime, Sharma2020MoistureProperties}). This can be divided into different regions, such as the flow region, the entanglement region, the transition region, and the glassy region \cite{Baurngaertel1992TheMelts}. The arrangement, length and entanglement of the polymer chains at the molecular level explains the different stiffness during the relaxation spectrum. A shorter chain relaxes faster than a longer chain and therefore does not contribute to the stiffness of the material after its relaxation \cite{Berry2006TheSolutions}.
To incorporate this behavior into our model, we will use an arbitrary, unknown number instead of a single Maxwell element. This leads us to the generalized Maxwell model, also known as the \emph{Maxwell-Wiechert model} (compare Figure \ref{abb:rheologicalModel}).\\

The number of Maxwell elements can be expanded to any number $n$ with the relaxation time $\tau_j$ in each of the dampers and the stiffness of the spring $\mu_j$ in each of the Maxwell elements. The stiffness of the individual spring is denoted by $\mu$. If all Maxwell elements relax to zero stress value, the equilibrium position is reached, and the single spring with stiffness $\mu$ represents the basic stiffness of the material, which ensures that the material does not behave like a liquid.
The deformation of the material represented by the strain value $\varepsilon$ is divided into an elastic component $\varepsilon_j^{\,e}$ and an inelastic component $\varepsilon_j^{\,i}$ in each of the Maxwell elements. The elastic component corresponds to the strain of the spring and the inelastic component to the strain of the damper. Since the strain of the damper depends on its relaxation time, an evolution equation based on the entropy principle is used to model the change of the inelastic strain with time  \cite{Johlitz2013Chemo-thermomechanicalMechanics, reese1998}. For small deformations, the evolution equation is given as
\begin{equation} 
\dot{\varepsilon}_j^{i}(t) = \frac{\varepsilon(t) - \varepsilon^{\,i}_j(t)}{\tau_j/2},
\label{eqn:evolution} 
\end{equation} 
where $\dot{\varepsilon}_j^{i}$ represents the time derivative. 
The total stress generated in the system is then given by the sum of the stresses induced in each of the springs. Assuming linear elastic behavior of the springs, this results in
\begin{equation}
\sigma(t) = \mu \, \varepsilon(t) + \sum\limits_{j=1}^{n}\, \mu_j (\varepsilon(t) - \varepsilon_j^{\,i}(t)).
\label{eqn:stressComputation}
\end{equation}

The total stress \eqref{eqn:stressComputation} is the sum of the stress of the individual spring, as well as the various Maxwell elements. The latter depends on the inelastic components of the strain, which are determined by the evolution equation \eqref{eqn:evolution}.

Apart from the mechanics \cite{reese1998}, the evolution equation \eqref{eqn:evolution} is also used in other applications such as Magnetic Particle Imaging \cite{croft2012relaxation} to model relaxation. Its solution can be formulated analytically as
\begin{align*}
\varepsilon_j^i(t)=\int \limits_0^t \varepsilon(\tilde{t}) \frac{2}{\tau_j} \exp\left( -2 \frac{t-\tilde{t}}{\tau_j} \right) d\tilde{t}.
 \end{align*}
If we use this solution, we can also represent the total stress analytically (cf. \cite{Rothermel2022}).
The stress $\sigma_j$ of the $j$-th Maxwell element with $j>1$ and the corresponding stiffness $\mu_j$ and relaxation time $\tau_j$ at time $t\in [0,T]$ is then represented as 
\begin{align}
\sigma_j(t)= \left\{\begin{array}{ll}  \frac{\mu_j \tau_j \dot{\varepsilon}}{2} \left( 1- \exp \left( -\frac{2}{\tau_j} t \right) \right) , &0 \leq t\leq \frac{\bar{\varepsilon}}{\dot{\varepsilon}} \\\
       - \frac{\mu_j \tau_j \dot{\varepsilon}}{2}  \left( 1- \exp \left( \frac{2 \bar{\varepsilon} }{\tau_j \dot{\varepsilon} } \right) \right) \exp \left( -\frac{2}{\tau_j} t \right), & \frac{\bar{\varepsilon}}{\dot{\varepsilon}} <t \leq T. \end{array}\right. 
        \label{eqn:stress_j}
\end{align}
The stress of the single spring, denoted by $\sigma_{0}$, can be specified with the corresponding stiffness $\mu$ as
\begin{align}
    \sigma_0(t)= \left\{\begin{array}{ll} \mu \dot{\varepsilon} t , &0 \leq t\leq \frac{\bar{\varepsilon}}{\dot{\varepsilon}} \\[1ex]
       \mu \bar{\varepsilon}, & \frac{\bar{\varepsilon}}{\dot{\varepsilon}} <t \leq T. \end{array}\right. 
       \label{eq:sigma0}
\end{align}

According to these calculations, the total stress can be written as the sum of the stresses of the single spring and the Maxwell elements, i.e.
\begin{align}
    \sigma (t)= \sum \limits_{j=0}^n \sigma_j(t) = \left\{\begin{array}{ll} \mu \dot{\varepsilon} t + \sum \limits_{j=1}^n \frac{\mu_j \tau_j \dot{\varepsilon}}{2} \left( 1- \exp\left( -\frac{2}{\tau_j} t \right) \right) , &0 \leq t\leq \frac{\bar{\varepsilon}}{\dot{\varepsilon}} \\
      \mu \bar{\varepsilon} - \sum \limits_{j=1}^n \frac{\mu_j \tau_j \dot{\varepsilon}}{2}  \left( 1- \exp \left( \frac{2 \bar{\varepsilon} }{\tau_j \dot{\varepsilon} } \right) \right) \exp \left( -\frac{2}{\tau_j} t \right) , & \frac{\bar{\varepsilon}}{\dot{\varepsilon}} <t \leq T. \end{array}\right. 
      \label{eq: sigmaTotal}
\end{align}

Thus, the inverse problem consists of identifying the material parameters from measured stress using equation \eqref{eqn:stressComputation}. The material parameters include the stiffness of the individual spring $\mu$, as well as the stiffnesses $\mu_1,\dots, \mu_n$ and relaxation times $\tau_1,\dots, \tau_n$ of the Maxwell elements. We emphasize that in our problem setting the number of Maxwell elements $n$ is assumed to be unknown and has to be computed as part of the solution, what, on the other hand, affects the number of material parameters to be determined. In this view the investigated inverse problem goes beyond existing approaches. The strain $\varepsilon$ is known and the inelastic strain $\varepsilon_j^i$ must be determined for all Maxwell elements $i=1,\dots,n$ from equation \eqref{eqn:evolution}.

%%%%%%%%%%%%%%%%%%%%%%%%%%%%%%%%%%%%%%%%%%%%%%%%%%

\subsection{The forward operator} \label{sec: Vorwärtsoperator}

%%%%%%%%%%%%%%%%%%%%%%%%%%%%%%%%%%%%%%%%%%%%%%%%%%%

The forward mapping associated with the inverse problem describes in section \ref{sec_visco} is given by
\begin{align*}
F: \mathcal{D}(F) \subset  \mathbb{N} \times \ell^2(\mathbb{N}) \to Y,\qquad  (n,x) \mapsto F(n,x):=\sigma,
\end{align*}
where the sequence $x=(\mu, \mu_1, \tau_1,\ldots,\mu_n, \tau_n, 0,0,\ldots)\in \l2$ contains the material parameters. Hence, $F$ maps $n\in \mathbb{N}$ and $x\in \l2 $ to the stress

\begin{align}
    \sigma (t)=      \left\{\begin{array}{ll} \mu \dot{\varepsilon} t + \sum \limits_{j=1}^n \frac{\mu_j \tau_j \dot{\varepsilon}}{2} \left( 1- \exp\left( -\frac{2}{\tau_j} t \right)  \right) , &0 \leq t\leq \frac{\bar{\varepsilon}}{\dot{\varepsilon}} \\
      \mu \bar{\varepsilon} - \sum \limits_{j=1}^n \frac{\mu_j \tau_j \dot{\varepsilon}}{2}  \left( 1- \exp \left( \frac{2 \bar{\varepsilon} }{\tau_j \dot{\varepsilon} }  \right)  \right) \exp \left( -\frac{2}{\tau_j} t \right) , &  \frac{\bar{\varepsilon}}{\dot{\varepsilon}} <t \leq T \end{array}\right.
      \label{eq:sigma}
\end{align}
with $\bar{\varepsilon},\dot{\varepsilon} > 0$.
Here, $n$ denotes the number of Maxwell elements. Since the number of material parameters depends on $n$, we cannot predict a priori how many parameters are to be determined. This is why we choose $x\in \l2$ allowing for a variable number of material parameters. Since there is only a finite number of parameters that are greater than zero, we set
\begin{align*}
    \mathcal{D}(F) = \lbrace (n,x) \in  I \times \ell^2(\mathbb{N}) : \quad 
    & I \subset \mathbb{N} \text{ finite}; \\
    & x:=\{x_m\}_{m\in \mathbb{N}},\ x_m \in \mathbb{R}_0^+,\ m \in \mathbb{N}; \\
    &x_m= 0,\ m> 2n+1; \\ 
    &x_{2i+1} \geq  \gamma,\ i=1,\dots,n, \ \gamma>0 \text{ arbitrarily small, but fixed}  \rbrace. 
\end{align*}

Here $I\subset \mathbb{N}$ denotes a physically reasonable finite subset of integers representing the possible number of Maxwell elements. 
%Since we are dealing with a nonlinear inverse problem, it is necessary to be in a neighborhood of the correct solution. This is especially true for the number of Maxwell elements $n$. For instance, if $n^*$ were the exact solution, $I:=\lbrace n^*-l, \dots, n^*+l \rbrace$ for $l\in \mathbb{N}$ would be an appropriate choice. \\
%Using the notation $x:=\{x_m\}_{m\in \mathbb{N}}$ we have
%\begin{align*}
%    x=\left(\mu, \mu_1, \tau_1,\dots,\mu_n, \tau_n, 0,0,\dots \right).
%\end{align*}
%implying $x_m=0$ for $m>2n+1$. 
The stiffnesses $\mu, \mu_1, \dots, \mu_n$ are non-negative.
From \eqref{eq:sigma} it follows that $\tau_i>0$ for $i=1,\dots,n$. So, there exists an (artificial) parameter $\gamma>0$ with $\tau_i \geq \gamma$ implying that $\mathcal{D}(F)$ is closed. This will be relevant in section \ref{sec:regularization_algorithm}. Without loss of generality we stick to the convention $\tau_1\leq\tau_2\leq\dots\leq\tau_n$.

%We will continue to use the rheological notation in the following and denote the material parameters by $\mu,\mu_j$ and $\tau_j$. The notation by the sequence elements $x_i$ will be taken up again in section \ref{sec:mathTheory}, where we will deal with the mathematical theory behind our developed procedure.\ 

%We adhere to the convention that $\tau_1\leq\tau_2\leq\dots\leq\tau_n$ holds and the Maxwell elements are ordered accordingly. \\

For fixed $n$ we denote by $F_n:\mathcal{D}(F_n) \subset \ell^2(\mathbb{N}) \to L^2([0,T])$ the mapping $F_n(x):=F(n,x)$.
Then,
\begin{align*}
    \mathcal{D}(F_n) = \lbrace &x:= \{x_m\}_{m\in \mathbb{N}}\in  \ell^2(\mathbb{N}) : \ x_m \in \mathbb{R}_0^+,\ m \in \mathbb{N}; \\
    &x_m= 0,\  m> 2n+1; \\ 
    &x_{2i+1} \geq  \gamma,\ i=1,\dots,n,\ \gamma>0 \text{ arbitrarily small, but fixed}  \rbrace .
\end{align*}

%%%%%%%%%%%%%%%%%%%%%%%%%%%%%%%%%%%%%%%%%%%%%%%%%%

\subsection{The inverse problem} \label{sec: inverse problem}

%%%%%%%%%%%%%%%%%%%%%%%%%%%%%%%%%%%%%%%%%%%%%%%%%%

The inverse problem is formulated as computing $n$ and $x$ as solution of
\begin{align}
\label{eq:unsere_operatorgleichung}
    F(n,x)=\sigma^\delta
\end{align}
from given, maybe noise-contaminated, stress measurements $\sigma^\delta$ with
$\|\sigma^\delta-\sigma\|\leq \delta$.
Note, that $F$ depends on the discrete variable $n$ as well as on the sequence $x$, where $n$ especially represents the number of material parameters to be determined corresponding to the non-zero elements in $x \in \l2$. Thus the number of material parameters depends on parts of the solution. This is an unusual situation in the field of inverse problems and differs from the theory in classical textbooks \cite{Engl1996RegularizationProblems, Louis1989InverseProbleme, Rieder2003KeineProblemen, Kaltenbacher2008IterativeProblems}. We develop an iterative solver relying on statistical inversion theory leading to the minimization of a Tikhonov functional 
\begin{align}
        T_\alpha(n,x):= \frac{1}{2}\lVert F(n,x) - \sigma^\delta \rVert^2 +\alpha \Omega(n,x)  
\end{align} 
where $\alpha>0$ is a regularization parameter and $\Omega$ represents an appropriate penalty term.

%%%%%%%%%%%%%%%%%%%%%%%%%%%%%%%

\section{A Bayes inversion algorithm} \label{chap: Bayes Inversion}
%%%%%%%%%%%%%%%%%%%%%%%%%%%%%%%%

%mit Beschreibung des Bayes Ansatzes in verbindung mit Tikhonov sowie der daraus resultierenden Iterationsmethode

%%%%%%%%%%%%%%%%%%%%%%%%%%%%%%%

In this section we develop a regularization method for solving \eqref{eq:unsere_operatorgleichung} relying on Bayes inversion and minimization of a resulting Tikhonov functional with a binomial prior. We give a brief introduction to statistical inversion theory, outline the application to \eqref{eq:unsere_operatorgleichung}, study the influence of the success probability to the prior and finally present our resulting iteration scheme.

%%%%%%%%%%%%%%%%%%%%%%%%%%%%%%%%%%%%%%%%%%%%%%%%%%

\subsection{Introduction to statistical inversion theory}

%%%%%%%%%%%%%%%%%%%%%%%%%%%%%%%%%%%%%%%%%%%%%%%%%%

The following introduction to statistical inversion theory is based on the textbooks \cite{Kaipio2005} and \cite{Dashti2017}. 

The goal of statistical inversion theory is to extract information about unobservable variables and to evaluate the uncertainty of these variables, based on all available knowledge about the measurement process and the information and models of the unknowns that are available before the measurement process.
The statistical inversion approach is based on the following principles:
\begin{enumerate}
    \item All variables, that are included in the model, are modeled as random variables.
    \item The randomness describes the degree of information about their
realizations.
\item The degree of information about these values is encoded in the probability distributions.
\item The result of the modeling is the a posteriori probability distribution.
\item The solution to the inverse problem is obtained by an estimator for this probability distribution, such as the maximum a posteriori (MAP) estimator.
\end{enumerate}

Let us consider a model 
\begin{align}
    f(x,e)=y.
    \label{eq: Bayes original model}
\end{align}
Here, $y\in \mathbb{R}^m$ denotes the measured quantity and $x\in \mathbb{R}^n$ is the unknown solution to be recovered. The mapping $f:\mathbb{R}^n \times \mathbb{R}^k \to \mathbb{R}^m$ represents the mathematical model, which may contain uncertainties and further unknown parameters. Moreover, let the measurement $y$ be contaminated by noise. The vector $e \in \mathbb{R}^k$ represents these additional unknowns as well as the noise.
As mentioned before, all variables of the model are modeled as random variables, i.e., \eqref{eq: Bayes original model} is to be seen as a realization of
\begin{align*}
    f(X,E)=Y
\end{align*}
with corresponding random variables $X$, $E$ and $Y$. Since the random variables are stochastically dependent on each other in this model, so are their random distributions. By $Y$ we denote the measurement random variable, where $Y=y_{\mathrm{meas}}$ describes its realization. The variable $E$ models the noise and $X$ the unknown parameters. We use the standard notation of writing a random variable in uppercase letters and its realization in lowercase.
A priori available information about $X$ is encoded by the probability density $\rho_0(x)=\rho_0(X=x)$, which is called the \emph{prior (distribution)}. We denote the probability density of the measurement by $\rho(y)$.
We assume that after analyzing the measurement and all other information available about the variables, the joint probability density of $X$ and $Y$ is given and denoted as $\rho(x,y)$. Then, the prior $\rho_0$ computes as
\begin{align*}
    \rho_0(x)=\int_{\mathbb{R}^m} \rho(x,y) dy.
\end{align*}
The conditional probability density of measurement $Y$ given $X=x$ computes as 
\begin{align*}
    \rho(y|x)=\frac{\rho(x,y)}{\rho_0(x)}
\end{align*}
for $\rho_0(x) > 0$, where $\rho(y|x)$ is called the \emph{likelihood function}. In the same way we determine the conditional probability 
\begin{align*}
    \rho(x|y_{\mathrm{meas}})= \frac{\rho(x,y_{\mathrm{meas}})}{\rho(y_{\mathrm{meas}})} 
\end{align*}
for $\rho(y_{mess})=\int_{\mathbb{R}^n} \rho(x,y_{mess})dx > 0$ which is called the \emph{a posteriori distribution}.
%Thus, in Bayesian statistics, an inverse problem can be expressed by the following task: for given data $Y=y_{measure}$ find the conditional probability $\rho(x|y_{measure})$ of the variable $X$. \\
A key ingredient of the stochastic approach to inverse problems is Bayes' Theorem.

\begin{theorem}[Theorem of Bayes]
    \label{theorem: Bayes theorem}
Let $X$ be a random variable with values in $\mathbb{R}^n$ and prior $\rho_0(x)$. In addition, let $Y$ be a random variable with values in $\mathbb{R}^k$ and realization $Y=y_{mess}$ with $\rho(y_{mess})>0$. Then, the a posteriori probability distribution of $X$ conditional on the data $y_{mess}$ is given by 
    \begin{align}
        \rho_{post}(x)=\rho(x|y_{mess})= \frac{\rho_0(x)\rho(y_{mess}|x)}{\rho(y_{mess})}.
        \label{eq: Bayes theorem}
    \end{align}
\end{theorem}

In the remainder of this article we simply write $y$ instead of $y_{mess}$ for simplicity.

Summarizing, there are three steps to solve an inverse problem using Bayes' formula \eqref{eq: Bayes theorem}.
First one has to find a probability density $\rho_0$ that contains all the relevant information about the variable $X$ that is available prior to the measurement process.
The second step is to choose the likelihood function $\rho(y|x)$ to describe the relationship between the measurement and the unknown. This contains the forward model, which is also used in classical inversion theory as well as any information about measurement or model error. 
Mostly, noise is modeled additively and independently of $X$. This is the approach that we will also apply to \eqref{eq:unsere_operatorgleichung}. We obtain the stochastic model
\begin{align}
    Y=f(X)+E,
    \label{eq: statistical model with noise}
\end{align}
with random variable $X$ with values in $\mathbb{R}^n$ and random variables $Y,E$ with values in $\mathbb{R}^m$, where $X$ and $E$ are stochastically independent. 
Assuming that the probability density of the noise is given by $\rho_{\text{noise}}(e)$ and using model \eqref{eq: statistical model with noise}, we conclude that $Y$ is identically distributed as $E$ conditioned on $X=x$, and the probability density of $E$ is shifted by $f(x)$. That means that the likelihood function can be expressed by 
\begin{align}
    \rho(y|x)=\rho_{\text{noise}}(y-f(x)).
    \label{eq: likelihood function general}
\end{align}
Thus, if we use \eqref{eq: Bayes theorem}, then we get the following expression for the a posteriori distribution:
\begin{align*}
    \rho(x|y) \propto \rho_0(x) \rho_{\text{noise}}(y-f(x)).
\end{align*}
In case that $X$ and $E$ are not stochastically independent or other noise models are considered than additive ones, we refer to \cite[p.56 ff]{Kaipio2005}.
Lastly, the a posteriori distribution has to be evaluated in a suitable way. A general distinction is made between \emph{point estimates} and \emph{interval estimates}. An interval estimate addresses the task to predict in which interval the values of the unknown is located at a certain probability given the prior $\rho_0$ and data $y$. A point estimate, on the other hand, determines a realization of $X$ by computing the expected value or maximizing the a posteriori density given the prior $\rho_0$ and data $y$. The last approach leads to the \emph{maximum a posteriori (MAP) estimator}. It is defined as
\begin{align*}
    x_{\text{MAP}}=\arg \max \limits_{x\in \mathbb{R}^n} \rho(x|y)
\end{align*} 
Note that, even if the maximizer exists, it needs not to be unique.

%%%%%%%%%%%%%%%%%%%%%%%%%%%%%%%%%%%%%%%%%%%%%%%%%%

\subsection{Application to the inverse problem} \label{sec: Bayes for us}

%%%%%%%%%%%%%%%%%%%%%%%%%%%%%%%%%%%%%%%%%%%%%%%%%%

In this section we use Bayes inversion theory to develop a regularization method for \eqref{eq:unsere_operatorgleichung}. In the following we at first focus on the reconstruction of $n$, the number of Maxwell elements.\\

Let $\Sigma$, $N$, $X$ and $E$ be random variables, we write $\sigma$ for the realization $\Sigma=\sigma$. 
Thus, the inverse problem \eqref{eq:unsere_operatorgleichung} can then be written as
 \begin{align*}
      \Sigma=F(N,X,E)
 \end{align*}
with $\Sigma$ the random variable for the measured data. The random variable $N$ attains integer values and describes the number of Maxwell elements, $X$ the material parameters and $E$ the noise.
We determine the MAP estimator using Bayes' Theorem \ref{theorem: Bayes theorem} and obtain
\begin{align*}
        n_{\text{MAP}}=\arg \max_{n\in \mathbb{N}} \rho(n|\sigma)\
        =\arg \max_{n\in \mathbb{N}}\frac{\rho_0(n)\rho(\sigma|n)}{\rho(\sigma)}.
    \end{align*}

The probability density $\rho(\sigma)$ of $\Sigma$ is independent of $n$ and thus irrelevant for computing the maximum. Thus, we obtain the maximization problem
\begin{align*}
        n_{\text{MAP}}=\arg \max_{n\in \mathbb{N}}\rho_0(n)\rho(\sigma|n).
    \end{align*}
We convert this into a minimization problem by applying $-\log(\cdot)$ leading to
\begin{align*}
     n_{\text{MAP}}= \arg \min_{n\in \mathbb{N}} \lbrace - \log (\rho(\sigma|n)) - \log (\rho_0(n)) \rbrace .
\end{align*}
We introduce the notation $\Phi(n|\sigma):=- \log (\rho(n|\sigma)) $, $\phi(\sigma|n):=- \log (\rho(\sigma|n))$ and $\phi_0(n):=- \log (\rho_0(n))$ turning the minimization problem to 
\begin{align*}
     n_{\text{MAP}}= \arg \min_{n\in \mathbb{N}} \Phi(n|\sigma) = \arg \min_{n\in \mathbb{N}} \lbrace \phi(\sigma|n) + \phi_0(n) \rbrace .
\end{align*}
First, we address the construction of the likelihood function $\rho(\sigma|n)$.
In our model we consider additive noise 
 \begin{align*}
      \Sigma=F(N,X)+E 
 \end{align*}
and assume that $N$ and $E$ as well as $X$ and $E$ are stochastically independent. Let $\rho_{\text{noise}}(e)$ be the probability density of the noise. Then, from \eqref{eq: likelihood function general}, we conclude 
\begin{align*}
    \rho(\sigma|n)=\rho_{\text{noise}}(\sigma-F(n,x)).
\end{align*}
In our application we suppose $E \sim \mathcal{N}(0,a^2)$, that is, 
\[  
\rho_{\text{noise}}(e)=\frac{1}{\sqrt{2\pi a^2}}\exp\left( \frac{-\lVert e \rVert ^2}{2a^2} \right).
\]
Applying $-\log(\cdot)$ leads to 
\begin{align*}
    -\log \left( \rho_{\text{noise}}(e) \right)=-\log\left( \frac{1}{\sqrt{2\pi a^2}} \right) + \frac{1}{2a^2} \lVert e \rVert ^2.
\end{align*}
Thus, the likelihood function reads as
\begin{align*}
    \phi(\sigma|n)=\phi_{\text{noise}}(\sigma-F(n,x))=-\log \left(\rho_{\text{noise}}(\sigma-F(n,x)) \right) \propto \lVert \sigma - F(n,x) \rVert^2.
\end{align*}
Minimizing $\phi(\sigma|n)$ for $n$, thus, is equivalent to minimizing $\lVert \sigma - F(n,x) \rVert^2$.

As prior we use the binomial distribution. It describes the number of successes in a series of independent trials that have two possible outcomes. Let $M$ be the number of trials and $0<q<1$ the probability of success in a single trial. Then, 
\begin{align*}
    B(n|q,M)= {M\choose n}q^n (1-q)^{M-n}
\end{align*}
is the probability of achieving $n$ successes. The usage of a binomial distribution for $n$ is motivated by the common strategy of subdividing the interval of relaxation times into subintervals of different decades and assuming that each decade contains a Maxwell element, see \cite{Scheffer2013Optimisationmaterials}. In our approach, the probability of every of the $M$ intervals containing an element is given by $q$ yielding more flexibility in the model. So, the approach can be described as guessing the number of Maxwell elements in a clever way and subsequently determining the material parameter by minimizing $\phi(\sigma | n)$ for $x$. To do this, let $I\subset \mathbb{N}$ a finite set of integers and $M:=\max I$ the maximum possible number of Maxwell elements. Then,
the prior is given as $\rho_0(n) = B(n|q,M)$ and $\phi_0(n)=-\log (\rho_0(n))$. 
Thus, minimizing 
\begin{align*}
     n_{\text{MAP}}= \arg \min_{n\in \mathbb{N}} \left\{ \lVert \sigma^\delta - F(n,x) \rVert^2 - \log \left({M\choose n}q^n (1-q)^{M-n}\right) \right\} 
\end{align*}
with $0<q<1$ results in finding the maximum a posteriori distribution $n_{\text{MAP}}$. To evaluate this functional, $x$ must be calculated for each considered $n$. Section \ref{sec: Bayes Algo} explains this step in more detail.
%If we look at this minimization from the point of view of classical inversion theory, we have found a penalty term using Bayes' theorem, which gives us a better way to estimate the number of Maxwell elements. 
To control the influence of the penalty term, we use a regularization parameter $\alpha >0$ and change the minimization finally to  
\begin{align}
\label{eq: min n_MAP}
     n_{\text{MAP}}= \arg \min_{n\in \mathbb{N}} \left\{ \lVert \sigma^\delta - F(n,x) \rVert^2 - \alpha \log \left({M\choose n}q^n (1-q)^{M-n}\right) \right\}. 
\end{align}

%%%%%%%%%%%%%%%%%%%%%%%%%%%%%%%%%%%%
\subsubsection{Impact of $q$ to the minimization} \label{sec: effect q prior}
%%%%%%%%%%%%%%%%%%%%%%%%%%%%%%%%%%%%

\begin{figure}
     \centering
    \subfloat{\def\svgwidth{0.5\textwidth}
   \hspace{-0.15cm}  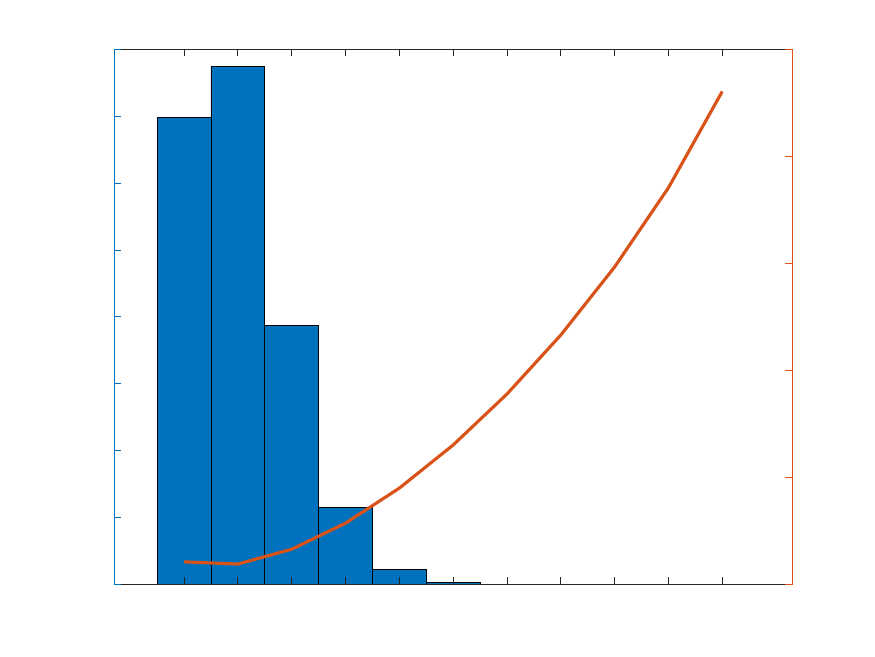}
 \subfloat{\def\svgwidth{0.5\textwidth}
   \hspace{0.15cm}  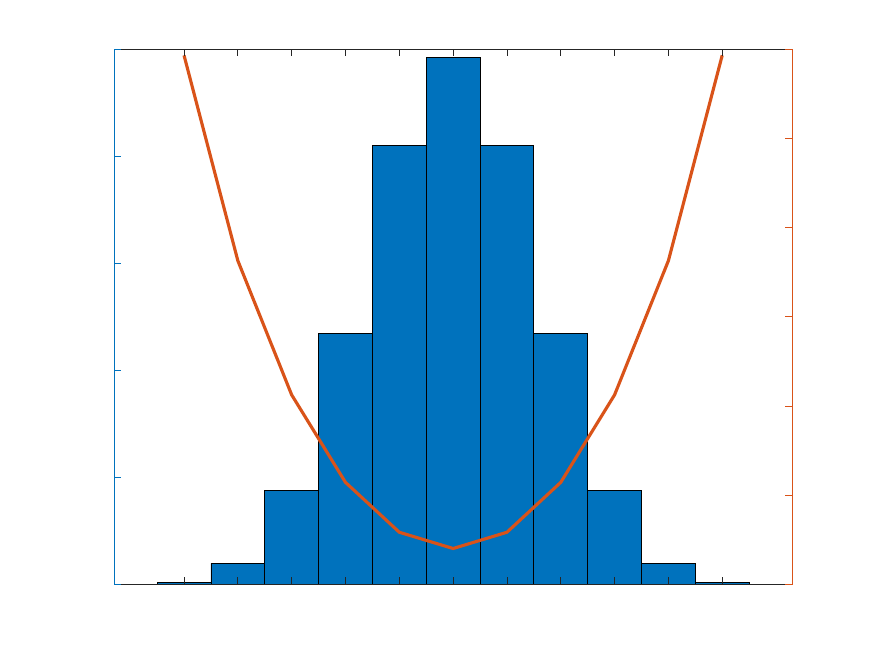}
   \caption{The binomial distribution $\rho_0(n)$ and $\phi_0(n)$ for $q=0.1$ (left picture) and $q=0.5$ (right picture)}
        \label{abb: Binomial q0.1 and q0.5}
\end{figure}

The probability density of the prior $\rho_0$ depends on the choice of success probability $q \in (0,1)$. We briefly investigate, how different values of $q$ affect the prior $\rho_0$ and, thus, $\phi_0$. 
For this purpose, we assume a maximum number of $M=10$ Maxwell elements.
On the left side of Figure \ref{abb: Binomial q0.1 and q0.5} one can see the binomial distribution $\rho_0(n)=B(n|q=0.1,M=10)$ (blue bars), and $\phi_0(n)=-\log(B(n|0.1,10))$ (red line) for $n=0,\dots, M$. The right picture shows according plots for $q=0.5$.

Since $\phi_0(n)$ increases with $n$, a large number of Maxwell elements is penalized in \eqref{eq: min n_MAP}. Figure \ref{abb: Binomial q0.1 and q0.5} shows the priors $\rho_0 (n)$ and the according penalties $\phi_0 (n)$ for different values of $q=0.1$ and $q=0.5$. Obviously small success probabilities $q$ penalize a big number of Maxwell elements significantly, whereas larger $q$ favor larger $n$. So, in order to get a good fit of the measured data by only a few number of elements it is recommendable to keep $q$ small. The numerical results in section \ref{sec:num_effect_q} will confirm that. In this view $q$ acts as an additional regularization parameter enforcing sparsity with respect to $x$ for small $q$.

%For instance, in the left picture of Figure \ref{abb: Binomial q0.1 and q0.5} $\rho_0(1)=0.3874$ is the highest probability and $\rho_0(10)\approx 0$ is the lowest. This is reflected in the minimization \eqref{eq: min n_MAP} by the values $\phi_0(1)\approx 0.95$ and $\phi_0(10) \approx 23.03$. Denoting, e.g., by $x_{1}$ the material parameters corresponding to a single Maxwell element and by $x_{10}$ the parameters corresponding to ten Maxwell elements the equivalence 
%\begin{align*}
%    \lVert \sigma^\delta - F(10,x_{1})\rVert^2 + \phi_0(1)  &< \lVert \sigma^\delta - F(1,x_{10}) \rVert^2 + \phi_0(10) \\
%    \Leftrightarrow  \qquad \lVert \sigma^\delta - F(10,x_{1})\rVert^2   &< \lVert \sigma^\delta - F(1,x_{10}) \rVert^2   -22.08
%\end{align*}
%shows that $n=10$ 
%must hold for the reconstruction with ten Maxwell elements to be equivalent to the solution of the MAP estimator for the given data. \\
%
%    
%On the right-hand side of Figure \ref{abb: Binomial q0.1 and q0.5} we see the prior $\rho_0$ for $q=0.5$. We see that for $q=0.5$ the mean number of Maxwell elements in this case $n=5$ is preferred. For higher values of $q$, correspondingly higher numbers of Maxwell elements are preferred. In most cases, it makes sense to choose $q$ smaller and thus favor a lower number of Maxwell elements in order to obtain the best possible description of the data by as few material parameters as feasible.

%%%%%%%%%%%%%%%%%%%%%%%%%%%%%%%%%%%%%%%%%%%%%%%%%%%%%
\subsubsection{The Bayes algorithm} \label{sec: Bayes Algo}
%%%%%%%%%%%%%%%%%%%%%%%%%%%%%%%%%%%%%%%%%%%%%%%%%%%%%

In this section we develop an iterative scheme for the minimization of 
\begin{align*}
    \min_{n\in \mathbb{N}} \Phi(n|\sigma)= \min_{n\in \mathbb{N}} \left\{ \lVert  F(n,x) - \sigma^\delta  \rVert^2 - \alpha\log \left({M\choose n}q^n (1-q)^{M-n}\right) \right\}.
\end{align*}

Let us define by $R:I \subset \mathbb{N} \to \mathcal{D}(F_n) \subset \l2$, $n \mapsto R(n)=x$ a minimization procedure that determines a minimizer $x$ of $\Phi(n|\sigma)$ for fixed $n$. This $x=R(n)$ is computed by minimizing a \emph{Tikhonov functional} of the form 
\begin{align}
\label{eq: Tikh in Bayes algo section}
       \min_{x \in \mathcal{D}(F_n)} T_{\beta,n}(x):= \min_{x \in \mathcal{D}(F_n)} \left\{ \frac{1}{2} \lVert F(n,x) - \sigma^\delta\rVert^2 +\alpha \Omega(x) \right\}
\end{align} 
with a penalty term $\Omega : \ell^2 (\mathbb{N})\to \mathbb{R}$, that serves as an additional stabilizing, and a corresponding regularization parameter $\alpha > 0$. We use an iterative solution method with several initial values for the stiffnesses and relaxation times 
$(\mu^{(0)}, \mu_1^{(0)}, \tau_1^{(0)}, \dots, \mu_n^{(0)}, \tau_n^{(0)},0,0 \dots )$. The initial values are distributed over the possible range of parameter values. The use of multiple initial values is necessary because of the nonlinearity of the problem causing the existence of multiple local minima. 
For each of these initial values, the Tikhonov functional \eqref{eq: Tikh in Bayes algo section} is minimized using the functions \textit{lsqnonlin} and \textit{MultiStart} from the MATLAB Optimization Toolbox \cite{Mathworks2016}. The function \textit{lsqnonlin} is a subspace trust-region method and is based on the interior-reflective Newton method described in \cite{Coleman1994, Coleman1996}.
As a stopping condition we use Morozov's discrepancy principle \cite{Morozov1984}. To this end we choose $\zeta > 1$ and stop the minimization as soon as 
$\lVert F_n(x) - \sigma^\delta \rVert \leq \zeta \delta$
with noise level $\delta > 0$. Finally we set $R(n):=x^\beta_n$ with the resulting minimizer $x^\alpha_n$ of $T_{\alpha,n}$.

The algorithm needs the following input values:
\begin{itemize}
\item the success probability $0<q<1$, where $q\to 1$ favors a large number of Maxwell elements $n$
\item measured data $\sigma^\delta$
\item $I:= \lbrace n_1,\dots,n_m \rbrace \subset \mathbb{N} $ with $n_1<\dots < n_m$, where $n_m=\max I =:M$ is the pre-selected maximum number of elements
\end{itemize}

As initialization of the algorithm we choose $k=1$ and $n_k:=n_1$, which corresponds to the smallest possible number of Maxwell elements. Subsequently, we compute $x_{1}^\alpha=R(n_1)$ and calculate $\Phi(n_1|\sigma)$. In the next step, we want to determine $l\in \lbrace 1, \dots, m-k \rbrace$ satisfying
\begin{align}
\label{eq: phi estimation}
    \Phi(n_{k+l}|\sigma) < \Phi(n_k|\sigma). 
\end{align}
This means that we have to compute $x_{k+l}^\alpha= R(n_{k+l})$ and $n_{k+l}$ from \eqref{eq: min n_MAP} for each $l$ to evaluate the expressions $ \Phi(n_{k+l}|\sigma)$. If an appropriate $l \in \lbrace 1, \dots, m-k \rbrace $ satisfying \eqref{eq: phi estimation} is not found, then we stop the iteration accepting $(n_k,x_k)$ as result. Otherwise we replace $k:=k+l$ and repeat the procedure to find $l$ with \eqref{eq: phi estimation} until the stopping criterion is achieved. Summarizing the algorithm reads as follows:

\begin{alg}\
\label{alg:Bayes}
Input: $0<q<1, M:=\max I:= \max \lbrace n_1,\dots,n_m \rbrace, \sigma \in L^2([0,T]), \zeta > 1$ (tolerance for $R$)
\begin{itemize}
  \item[1)] Let $k=1$ and $n_k:=n_1$, compute $x_k=x_1$ by $R(n_k)= x_k^\alpha$. 
  \item[2)] Determine $l\in \lbrace 1, \dots, m-k \rbrace$ with $\Phi(n_{k+l}|\sigma) < \Phi(n_k|\sigma)$ and set $k:=k+l$. For each $l$ calculate $x_{k+l}^\alpha= R(n_{k+l})$ to evaluate $\Phi(n_{k+l}|\sigma)$.
     \item[3)] If no such $l$ exists, then STOP. Otherwise proceed with step 2).
\end{itemize}
Output: $(n_k,x_k^\alpha)$
\end{alg}

Algorithm \ref{alg:Bayes} is an iteration scheme minimizing alternately $\Phi(n|\sigma)$ by \eqref{eq: min n_MAP} and $T_{\alpha,n}$ with respect to $x$. Since $I$ is finite, the existence of a minimizing $n_k$ is guaranteed and we obtain, that Algorithm \ref{alg:Bayes} iteratively computes a minimizer of the Tikhonov functional
\begin{align}\label{eq:Tikhonov_Bayes}
       \min_{(n,x) \in \mathcal{D}(F)} T_\alpha(n,x):= \min_{(n,x) \in \mathcal{D}(F)} \left\{ \frac{1}{2}\lVert F(n,x) - \sigma^\delta \rVert^2 +\alpha \left[ \Omega(x) -\log \left({M\choose n}q^n (1-q)^{M-n}\right) \right] \right\}
\end{align} 
with a (possible generic) regularization parameter $\alpha>0$.

%%%%%%%%%%%%%%%%%%%%%%%%%%%%%%%%%%%%%%%%%%%%%%%%%%%%%%%%%%%%%%%%%%%%%%%%%
%%%%%%%%%%%%%%%%%%%%%%%%%%%%%%%%%%%%%%%%%%%%%%%%%%%%%%%%%%%%%%%%%%%%%%%%%

\section{Convergence, stability and regularization of Algorithm \ref{alg:Bayes}}
\label{sec:regularization_algorithm}

This section is devoted to prove convergence, stability and regularization properties of Algorithm \ref{alg:Bayes}. To this end in Section \ref{sec:Tikhonov} we recapitulate well-known results for Tikhonov regularization of nonlinear inverse problems in Banach spaces. Section \ref{sec:semigroups} extends the well-established regularization theory to $\mathbb{N}\times \ell^2 (\mathbb{N})$ wich is necessary since it does not apply to discrete semigroups as $\mathbb{N}$. Using the new framework we are able to prove convergence, stability and regularization properties for algorithm \ref{alg:Bayes}.

%%%%%%%%%%%%%%%%%%%%%%%%%%%%%%%%%%

\subsection{Ill-posedness and regularization in Banach spaces}
\label{sec:Tikhonov}

%%%%%%%%%%%%%%%%%%%%%%%%%%%%%%%%%%

We briefly recall the established concepts of Toikhonov regularization in Banach spaces from \cite[Sections 3 and 4]{Schuster2012RegularizationSpaces} as they are important for our further investigations. For further insights we recommend \cite{Ito2015InverseAlgorithms, Hofmann2007AOperators, Seidman1989WellProblems}.

We consider a nonlinear operator equation 
\begin{align}
\label{eq:nonlinear operator equation}
    F(x)=y, \quad x \in \mathcal{D}(F) \subseteq X, \quad y \in F(\mathcal{D}(F))\subseteq Y, 
\end{align}
with domain $\mathcal{D}(F)$ and
$$F(\mathcal{D}(F))= \{ \tilde{y}\in Y: F(\tilde{x})=\tilde{y} \text{ for some } \tilde{x}\in \mathcal{D}(F) \}$$
the image of the nonlinear operator $F:\mathcal{D}(F) \subseteq X \to Y $. Here,
$X$ and $Y$ are Banach spaces. In any Banach space we consider strong (norm) convergence, denoted by $\to$ and the weak topology, where convergence is denoted by $\rightharpoonup$.

%%%%%%%%%%%%%%%%%%%%%%%%%%%%%%%%%%%%%%%%%%%%%%%%%%

\subsubsection{Local ill-posedness}

%%%%%%%%%%%%%%%%%%%%%%%%%%%%%%%%%%%%%%%%%%%%%%%%%%

%Corresponding to the local character of solutions in nonlinear equations, we also have to consider the ill-posedness locally, i.e. in the vicinity of a solution of the operator equation (compare \cite{Hofmann1994FactorsProblems, Hofmann2000Ill-posednessSpaces, Hofmann1998LocalSpaces}). 
%In contrast, for linear problems, ill-posedness could also be captured globally. 
%Further theory on ill-posedness in linear problems can be found in \cite{Engl1996RegularizationProblems,Louis1989InverseProbleme, Rieder2003KeineProblemen, Schuster2012RegularizationSpaces}. 

It is inherent to a locally ill-posed equation \eqref{eq:nonlinear operator equation} that a solution $x_0$ cannot be computed stably even if the perturbations on the right-hand side $y=F(x_0)$ are arbitrarily small. 

\begin{dfn}[Local ill-posedness]
\label{dfn:local ill-posedness}
A nonlinear operator equation \eqref{eq:nonlinear operator equation} is called \textit{local ill-posed} at $x_0 \in \mathcal{D}(F)$ satisfying $F(x_0)=y$ if for arbitrarily small radii $r>0$ there exists a sequence $\{x_n\}_{n=1}^{\infty} \subset \bar{B}_r(x_0) \cap \mathcal{D}(F)$ such that 
\begin{align*}
    F(x_n) \to F(x_0) \text{ in } Y, \quad \text{ but } x_n \nrightarrow x_0 \text{ in } X,\quad \text{ for } n \to \infty.  
\end{align*}
Otherwise, the equation is called \textit{locally well-posed} at $x_0$. 
\end{dfn}

Since in general only noisy data $y^\delta$ satisfying
\begin{align*}
    \lVert y^\delta - y \rVert \leq \delta 
\end{align*}
with noise level $\delta >0$ are available, local ill-posedness leads to instabilities in the solution process of \eqref{eq:nonlinear operator equation}.
A stable solution of nonlinear, ill-posed equations requires so-called \emph{reguarization methods}.
 
\begin{dfn}[Regularization method]
A mapping that assigns each pair $(y^\delta,\alpha) \in Y \times (0,\bar{\alpha})$ with $0< \bar{\alpha} \leq + \infty$ to a well-defined element $x^\delta_\alpha \in \mathcal{D}(F)$, is called \textit{regularization} (regularization method) for the nonlinear operator equation \eqref{eq:nonlinear operator equation}, if there is a suitable choice $\alpha=\alpha(y^\delta,\delta)$ of the regularization parameter such that for any sequence $\{ y_n \}_{n=1}^{\infty}\subset Y$ with $\lVert y_n-y \rVert \leq \delta_n$ and $\delta_n \to 0$ for $n \to \infty$, the corresponding regularized solutions $x_{\alpha(y_n,\delta_n)}^{\delta_n}$ converge with respect to a topology given on $X$ to the solution $x^+$ of the equation \eqref{eq:nonlinear operator equation}. If it is not unique, regularized solutions must converge to solutions of the equation \eqref{eq:nonlinear operator equation} which have desired properties, e.g., to $\bar{x}$-minimum norm solutions. In the case of non-uniqueness, different subsequences of regularized solutions may converge to different solutions of the operator equation, all having the same desired property. 
\end{dfn}

A comprehensive treatise of regularization methods can be found in the textbooks \cite{Engl1996RegularizationProblems, Louis1989InverseProbleme, Rieder2003KeineProblemen}. Nonlinear inverse problems are particularly subject of \cite{Kaltenbacher2008IterativeProblems}. 
In this article we focus on Tikhonov regularization, see, e.g. \cite{Engl1989ConvergenceProblems, Kaltenbacher2008IterativeProblems, Neubauer1992TikhonovScales, Scherzer1993TheProblems}.
Other well-known regularization methods are the Landweber method \cite{Hanke1995AProblems, Kaltenbacher1997SomeProblems, Kaltenbacher2009IterativeBanachSpacesProblems}, Newton type methods \cite{Kaltenbacher2009IterativeBanachSpacesProblems}, or the method of approximate inverse \cite{Louis1995ApproximateProblems}.

%%%%%%%%%%%%%%%%%%%%%%%%%%%%%%%%%%%%%%%%%%%%%%%%%%
\subsubsection{Tikhonov regularization}
\label{sec:Tikhonov-regularization}
%%%%%%%%%%%%%%%%%%%%%%%%%%%%%%%%%%%%%%%%%%%%%%%%%%

The Tikhonov regularization assigns a pair $(y^\delta,\alpha)$ to a (well-defined) minimizer $x_\alpha^\delta$ of the functional 
\begin{align}
    T_\alpha(x):= \frac{1}{p}\lVert F(x)-y^\delta \rVert^p + \alpha \Omega(x) 	 
	\quad  \text{ subject to } x\in \mathcal{D}(F)\subseteq X.
    \label{eq: Tikhonov}
\end{align}
By $\Omega:X \to [0,\infty]$, we denote a penalty term that guarantees stability and enforces desired properties to the solution $x_\alpha^\delta$. The first step is to prove that there always exists a well-defined minimizer to any pair $(y^\delta,\alpha)$. Subsequently we outline stability and convergence of Tikhonov regularization. To this end we define $\Omega$-minimizing solutions (c.f. \cite{Hofmann2007AOperators}).

\begin{dfn}[$\Omega$-minimizing solution]
\label{dfn: Omega min solution}
We say that $x^+\in \mathcal{D}(F) \subseteq X $ is a $\Omega$-minimizing solution of the operator equation $\eqref{eq:nonlinear operator equation}$ if
\begin{align*}
    F(x^+)=y \quad \text{ and } \quad \Omega(x^+)=\inf \{ \Omega(\tilde{x}) \, : \, \tilde{x}\in \mathcal{D}(F), \, F(\tilde{x})=y \}.
\end{align*}
\end{dfn} 

Next we formulate two assumption blocks \ref{Ann: 1st block} and \ref{Ann: 2nd block} which are essential for turning the minimization of Tikhonov functionals \eqref{eq: Tikhonov} to a regularization method.
The first block makes assumptions about the Banach spaces $X$ and $Y$, as well as the operator $F$ and its domain $\mathcal{D}(F)$. 

\begin{assump} \ \\
\label{Ann: 1st block}
\vspace{-0.7cm}
\begin{itemize} 
\item[(a)] $X$ and $Y$ are infinite dimensional reflexive Banach spaces.
\item[(b)] $\mathcal{D}(F)$ is a closed and convex subset of $X$.
\item[(c)] For $x_n \rightharpoonup x_0 $ in $X$, $x_n \in \mathcal{D}(F),\, n \in \mathbb{N},\, x_0 \in \mathcal{D}(F)$ follows $F(x_n) \rightharpoonup F(x_0)$ in $Y$, i.e., $F: \mathcal{D}(F) \subseteq X \to Y $ is weak-to-weak sequentially continuous.
\end{itemize}
\end{assump}

The next assumptions are associated with the penalty term $\Omega(x)$ of the Tikhonov functional $T_\alpha(x)$ and the $\Omega$-minimizing solution of the operator equation.

\begin{assump}\ \\
\label{Ann: 2nd block}
\vspace{-0.7cm}
\begin{itemize} 
\item[(a)] For the exponent $p$ in \eqref{eq: Tikhonov} holds $1<p<\infty$.
\item[(b)] $\Omega$ is a \emph{proper}, convex and lower semi-continuous functional, where proper means that the domain of $\Omega$ is nonempty. Furthermore 
\begin{align*}
    \mathcal{D}:= \mathcal{D}(F) \cap \mathcal{D}(\Omega) \neq \emptyset.
\end{align*}
\item[(c)] It is assumed that $\Omega$ is a stabilizing functional in the sense that the subsets 
\begin{align*}
\mathcal{M}_\Omega(c):=\{ x \in X \, : \, \Omega(x) \leq c \}
\end{align*}
in $X$ for all $c \geq 0$ are weakly sequentially pre-compact.
\item[(d)] There exists a $\Omega$-minimizing solution $x^+$ of the equation \eqref{eq:nonlinear operator equation}, which belongs to the so-called Bregman domain
\begin{align*}
\mathcal{D}_B(\Omega):=\lbrace x\in \mathcal{D} \subseteq X: \partial\Omega(x) \neq \emptyset \rbrace,
\end{align*}
 where $\partial \Omega (x) \subseteq X^*$ denotes the subdifferential of $\Omega$ at point $x$. 
\end{itemize}
\end{assump}

Penalty terms of the form 
\begin{align}
\label{eq: omega norm form}
   \Omega(x):= \frac{1}{q} \lVert x \rVert^q_X , \qquad 1 \leq q < \infty
\end{align}
are always stabilizing functionals, if $X$ is a reflexive Banach space as postulated by \ref{Ann: 1st block}(a).
The proof of this result can be found in \cite{Schuster2012RegularizationSpaces} or \cite[p.251]{Megginson1998AnTheory}. 

Assuming \ref{Ann: 1st block} and \ref{Ann: 2nd block} one can deduce the following main results. For the proof we refer again to \cite{Schuster2012RegularizationSpaces}. 

\begin{theorem}[Existence]
\label{theorem:existence}
For all $\alpha>0$ and $y^\delta\in Y$ there exists an element $x_\alpha^\delta \in \mathcal{D}(F)$ which minimizes the Tikhonov functional $T_\alpha(x)$ \eqref{eq: Tikhonov} in $\mathcal{D}(F)$. 
\end{theorem}

The next property we consider is stability with respect to the data $y^\delta$. 

\begin{theorem}[Stability]
\label{theorem: stability}
For all $\alpha>0$, the minimizers $x_\alpha^\delta$ of \eqref{eq: Tikhonov} are stable with respect to the data $y^\delta$. That is, for a sequence of data $\{ y_n \}$ converging to $y^\delta$ in the norm topology of $Y$, i.e., $$\lim_{n\to\infty}\lVert y_n-y^\delta \rVert_Y=0,$$ any associated sequence of minimizers $\{ x_n \}$ of the extremal problem 
\begin{align*}
    \frac{1}{p}\lVert F(x)-y_n \rVert^p + \alpha \Omega(x)\quad  \text{ subject to } x\in \mathcal{D}(F)\subseteq X
\end{align*}
has a subsequence $\{ x_{n_k} \}$ that converges in the weak topology of $X$, and the weak limit $\Tilde{x}$ of any such subsequence is a minimizer $x_\alpha^\delta$ of \eqref{eq: Tikhonov}. Moreover, for any such weakly convergent subsequence holds $$\lim \limits_{k\to \infty} \Omega(x_{n_k})=\Omega(x_\alpha^\delta).$$
\end{theorem}

As a next step we formulate conditions for convergence in case that the noise level $\delta \to 0$. 

\begin{theorem}[Convergence]
\label{theorem: convergence}
Let $\{ y_n:=y^{\delta_n} \}\subset Y$ be a sequence of perturbed data to the exact right-hand side $y\in F(\mathcal{D}(F))$ of equation \eqref{eq:nonlinear operator equation} and let $\lVert y-y_n\rVert \leq \delta_n$ for a sequence $\{ \delta_n >0 \}$ that converges monotonically to zero. Moreover, we consider a sequence $\{ \alpha_n >0 \}$ of regularization parameters and an associated sequence $\{ x_n:=x_{\alpha_n}^{\delta_n} \}$ of regularized solutions, which are minimizers of 
\begin{align*}
    \frac{1}{p} \lVert F(x)- y_n \rVert^p + \alpha_n\Omega(x) \text{ with } x\in \mathcal{D}(F)\subseteq X.
\end{align*}
Under the conditions 
\begin{align*}
    \limsup \limits_{n\to \infty} \Omega(x_n) \leq \Omega(x_0) \text{ for all } x_0 \in \{ x\in \mathcal{D}:= \mathcal{D}(F)\cap \mathcal{D}(\Omega) : F(x)=y \}
\end{align*}
and 
\begin{align*}
    \lim \limits_{n\to \infty} \lVert F(x_n)-y_n \rVert =0 
\end{align*}
the sequence $\{ x_n \}$ has a weakly convergent subsequence, where each weak limit of such a subsequence is a $\Omega$-minimizing solution $x^+\in \mathcal{D}$ of the operator equation. Moreover, if the $\Omega$-minimizing solution $x^+\in \mathcal{D}$ is unique, we obtain weak convergence $x_n \rightharpoonup x^+$ in $X$.
\end{theorem}

%%%%%%%%%%%%%%%%%%%%%%%%%%%%%%%%%%%%%%%%%%%%%%%%%%

\subsection{Convergence, stability and regularization in $\mathbb{N}\times \ell^2 (\mathbb{N})$} 
\label{sec:semigroups}

%%%%%%%%%%%%%%%%%%%%%%%%%%%%%%%%%%%%%%%%%%%%%%%%%%

In this section we apply the theoretical results from Section \ref{sec:Tikhonov-regularization} to algorithm \ref{alg:Bayes}. To this end we have to take into account that the Tikhonov functional \eqref{eq:Tikhonov_Bayes} is formulated in $\mathbb{N}\times \ell^2 (\mathbb{N})$ where $\mathbb{N}$ shows the structure of a discrete semigroup. This means that the existing theory has to be extended in a suitable way.

%%%%%%%%%%%%%%%%%%%%%%%%%%%%%%%%%%%%%%%%%%%%%%%%%%

\subsubsection{Local ill-posedness of \eqref{eq:unsere_operatorgleichung}}

We consider $\mathbb{N}\times \ell^2(\mathbb{N})$ with the product topologies $(\tau_d,\tau_0)$ and $(\tau_d,\tau_w)$. Here, $\tau_0$ represents the strong topology in $\ell^2(\mathbb{N})$ induced by the $\ell^2$-norm
\begin{align*}
    \lVert x \rVert_{\ell^2(\mathbb{N})}:= \left( \sum \limits_{n=1}^\infty x_n^2 \right)^{\frac{1}{2}} \quad \text{ for } x \in \ell^2(\mathbb{N})
\end{align*}
and $\tau_w$ denotes the weak topology in $\ell^2(\mathbb{N})$. 
Furthermore, $\tau_d$ denotes the discrete topology induced by the discrete metric
\begin{align*}
    d(x,y):= \left\{\begin{array}{ll} 0, & x=y \\
         1, & x\neq y \end{array}\right.  \quad \text{ for } x,y \in \mathbb{N}.
\end{align*}
A sequence $\{ (n_m,x_m) \}$ in $\mathbb{N}\times \ell^2(\mathbb{N})$ converges if and only if $\{ n_m \}$ converges in $\mathbb{N}$ and $\{ x_m \}$ in $\ell^2(\mathbb{N})$.
A sequence $\{ n_m \}\subset \mathbb{N}$ converges with respect to $\tau_d$ if and only if it is constant starting from a certain element, i.e., there exists an $m_0\in\mathbb{N}$ with $n_m = \mathrm{const}$ for all $m\geq m_0$.

By 
$$\mathcal{U}(n):=\{n-m, \dots, n+l\}$$
for some $m,l\in \mathbb{N}_0$ we denote a $\tau_d$-neighborhood of $n\in \mathbb{N}$ in $\mathbb{N}$.

\begin{dfn}
$F(n,x)=\sigma$ is called \emph{locally well-posed} in $(n^+,x^+)\in \mathcal{D}(F)$ if $F_n(x)=\sigma$ is locally well-posed for all $n\in\mathcal{U}(n^+)\cap I$ and all neighborhoods $\mathcal{U}(n^+)$ of $n^+$.
$F(n,x)=\sigma$ is called \emph{locally ill-posed} in $(n^+,x^+)\in \mathcal{D}(F)$ if there exists a neighborhood $\mathcal{U}(n^+)$ such that $F_n(x)=\sigma$ is locally ill-posed for some $n\in\mathcal{U}(n^+) $ according to definition \ref{dfn:local ill-posedness}.\\
\end{dfn}

\begin{theorem}
There exists a neighborhood $\mathcal{U}(n^+)$ of $n^+$ such that $F_n(x)=\sigma$ is locally ill-posed for some $n\in\mathcal{U}(n^+)\cap I$.\\  
\end{theorem}
\textit{Proof:}\\
Let $r>0$ and $x^+ = \left(\mu^+, \mu_1^+, \tau_1^+,\dots,\mu_{n^+}^+, \tau_{n^+}^+, 0,0,\dots \right) $ be fixed.
We choose $\mathcal{U}(n^+)$, such that $n=n^++1 \in\mathcal{U}(n^+)$, and $x^r_k= \left(\mu, \mu_1^+, \tau_1^+,\dots,\mu_{n^+}^+, \tau_{n^+}^+, \frac{r}{2k},\frac{r}{2},0,0\dots \right) $.
Then, $(n,x_k^r)\in \mathcal{D}(F)$ for all $k\in \mathbb{N}$.
Furthermore, it follows that
\begin{align*}
   \lVert x_k^r-x^+ \rVert_{\ell^2(\mathbb{N})}^2 = \frac{r^2}{4k^2} + \frac{r^2}{4} = \left( \frac{1}{4k^2} + \frac{1}{4} \right)r^2 <r^2 \text{ for all } k\in \mathbb{N}.
\end{align*}
This yields $\{x_k^r\}_{k\in \mathbb{N}} \subset B_r(x^+)$ and from 
\begin{align*}
   \lVert x_k^r-x^+ \rVert_{\ell^2(\mathbb{N})}^2 \to  \frac{r^2}{4} 
   \qquad \mbox{as } k\to\infty
\end{align*}
we deduce $x_k^r \nrightarrow x^+$ as $k \to \infty$. We compute
\begin{align*}
   \lVert F_n(x_k^r)&-F_{n^+}(x^+) \rVert_{L^2([0,T])}^2=\\
   \int \limits_0^{\bar{\varepsilon}/\dot{\varepsilon}} & 
   \left[ \mu^+ \dot{\varepsilon} t 
   + \sum \limits_{j=1}^n \frac{\mu_j^+ \tau_j^+ \dot{\varepsilon}}{2} \left( 1- \exp\left( -\frac{2}{\tau_j^+} t \right)  \right) 
   + \frac{\frac{r}{2k}  \frac{r}{2}\dot{\varepsilon}}{2} \left( 1- \exp\left( -\frac{2}{\frac{r}{2} } t \right)  \right) \right. \\
   & \left. - \mu^+ \dot{\varepsilon} t - \sum \limits_{j=1}^n \frac{\mu_j^+ \tau_j^+ \dot{\varepsilon}}{2} \left( 1- \exp\left( -\frac{2}{\tau_j^+} t \right)  \right)   \right]^2  dt \\
   + \int \limits_{\bar{\varepsilon}/\dot{\varepsilon}}^T 
   & \left[ \mu^+ \bar{\varepsilon} 
   - \sum \limits_{j=1}^n \frac{\mu_j^+ \tau_j^+ \dot{\varepsilon}}{2}  \left( 1- \exp \left( \frac{2 \bar{\varepsilon} }{\tau_j^+ \dot{\varepsilon} }  \right)  \right) \exp \left( -\frac{2}{\tau_j^+} t \right) 
   + \frac{\frac{r}{2k} \frac{r}{2} \dot{\varepsilon}}{2}  \left( 1- \exp \left( \frac{2 \bar{\varepsilon} }{\frac{r}{2} \dot{\varepsilon} }  \right)  \right) \exp \left( -\frac{2}{\frac{r}{2}} t \right) \right. \\
   & \left. - \mu^+ \bar{\varepsilon} - \sum \limits_{j=1}^n \frac{\mu_j^+ \tau_j^+ \dot{\varepsilon}}{2}  \left( 1- \exp \left( \frac{2 \bar{\varepsilon} }{\tau_j^+ \dot{\varepsilon} }  \right)  \right) \exp \left( -\frac{2}{\tau_j^+} t \right) \right]^2 dt \\
   = \int \limits_0^{\bar{\varepsilon}/\dot{\varepsilon}} &
   \left[ \frac{r^2 \dot{\varepsilon}}{8k} \left( 1- \exp\left( -\frac{4}{r} t \right)  \right) \right]^2 dt 
   + \int \limits_{\bar{\varepsilon}/\dot{\varepsilon}}^T  
   \left[ \frac{r^2 \dot{\varepsilon}}{8k}  \left( 1- \exp \left( \frac{4 \bar{\varepsilon} }{r \dot{\varepsilon} }  \right)  \right) \exp \left( -\frac{4}{r} t \right) \right]^2 dt \\
   =\frac{1}{k^2} & \left( \int \limits_0^{\bar{\varepsilon}/\dot{\varepsilon}}  \left[ \frac{  r^2\dot{\varepsilon}}{8} \left( 1- \exp\left( -\frac{4}{r} t \right)  \right) \right]^2 dt 
   + \int \limits_{\bar{\varepsilon}/\dot{\varepsilon}}^T  
   \left[ \frac{r^2\dot{\varepsilon}}{8}  \left( 1- \exp \left( \frac{4 \bar{\varepsilon} }{r \dot{\varepsilon} }  \right)  \right) \exp \left( -\frac{4}{r} t \right) \right]^2 dt \right). 
\end{align*}
Since $\left( 1- \exp\left( -\frac{4}{r} t \right)\right)^2$ and $\exp \left( -\frac{4}{r} t \right)$ are continuous, the integrals remain bounded and, thus, 
$$ \lVert F_n(x_k^r)-F_{n^+}(x^+) \rVert_{L^2([0,T])}^2 \to 0
\qquad \mbox{as } k\to\infty.$$ \qed

%%%%%%%%%%%%%%%%%%%%%%%%%%%%%%%%%%%%%%%%%%%%%%%%%%%%%%%%%%%%%%%%%%%%%

\subsubsection{Algorithm \ref{alg:Bayes} as a regularization method}

As outlined in section \ref{sec: Bayes Algo} we solve the inverse problem by minimizing \eqref{eq:Tikhonov_Bayes}, where $0<q<1$ is the success probability, $M$ is the maximum number of Maxwell elements, $\alpha$ is a regularization parameter, and $\Omega$ is a penalty term satisfying assumptions \ref{Ann: 2nd block}.
For simplicity, we define $$\logg(n):=-\log \left({M\choose n}q^n (1-q)^{M-n}\right).$$

%\textbf{Considering the subproblem}

For fixed $n\in I$ we consider the problem
\begin{align}
\label{eq:F_n equation}
    F_n(x)=\sigma, \quad x \in \mathcal{D}(F_n) \subseteq \ell^2(\mathbb{N}), \quad \sigma \in F_n(\mathcal{D}(F_n))\subseteq L^2([0,T]), 
\end{align}
with $\mathcal{D}(F_n)$ the domain of $F_n$ and $$F_n(\mathcal{D}(F_n))=\{\tilde{\sigma}\in L^2([0,T]): F_n(\tilde{x})=\tilde{\sigma} \text{ for a } \tilde{x}\in \mathcal{D}(F_n) \}$$ the image of $F_n$. The inverse problem \eqref{eq:F_n equation} is solved by minimizing the Tikhonov functional \eqref{eq: Tikh in Bayes algo section}.
We show that the assumptions \ref{Ann: 1st block} and \ref{Ann: 2nd block} are satisfied for \eqref{eq: Tikh in Bayes algo section} by proving that a regularized solution $x_\alpha^\delta$ exists, that it is stable with respect to the data, and that, given certain conditions, the regularized solution converges weakly.
After that, we use propositions \ref{theorem:existence}, \ref{theorem: stability}, and \ref{theorem: convergence} to also show these properties for the overall procedure \eqref{eq:Tikhonov_Bayes}.\\

Let us start by proving the assumptions \ref{Ann: 1st block} for $F_n:\mathcal{D}(F_n) \subset \ell^2(\mathbb{N}) \to L^2([0,T])$:
\begin{enumerate}[leftmargin=0.5cm]
    \item[(a)] $X=\ell^2(\mathbb{N}),\ Y= L^2([0,T])$ are infinite dimensional Hilbert spaces and, thus, reflexive. 
    \item[(b)] The domain 
\begin{align*}
    \mathcal{D}(F_n) = \lbrace &x:= \{x_m\}\in  \ell^2(\mathbb{N}) : \ x_m \in \mathbb{R}_0^+,\ m \in \mathbb{N}; \\
    &x_m= 0,\  m> 2n+1; \\ 
    &x_{2i+1} \geq  \gamma,\ i=1,\dots,n,\ \gamma>0 \text{arbitrarily small but fixed}   \rbrace 
\end{align*}
  is closed and convex. Let $\{ x^{(k)} \}_{k\in \mathbb{N}}$ be a sequence in $\mathcal{D}(F_n)\subset\ell^2(\mathbb{N})$. Hence, $x^{(k)}:=\{x_m \}^{(k)}_{m\in \mathbb{N}}$ for $k\in \mathbb{N}$.
    Moreover, let $x^{(k)} \to x^*$ with $x^* \in \ell^2(\mathbb{N})$ as $k\to +\infty$. We show that $x^* \in \mathcal{D}(F_n)$.
    From $x^{(k)} \to x^*$ as $k\to +\infty$ it follows that for all $\varepsilon >0$ there exists a $K\in \mathbb{N}$ such that for all $k\geq K$ it holds
    \begin{align*}
        \lVert x^{(k)}-x^* \rVert_{\l2}^2 < \varepsilon.
    \end{align*}
     Since $x^{(k)} \in \mathcal{D}(F_n)$, we have $x_m^{(k)}=0$ for all $k\in \mathbb{N}$ and $m>2n+1$, and thus 
     \begin{align*}
        \lVert x^{(k)}-x^* \rVert_{\l2}^2= \sum \limits_{m=1}^{2n+1} \lvert x_m^{(k)} -x_m^* \rvert^2 + \sum \limits_{m=2n+2}^{\infty} {(x_m^*)}^2< \varepsilon.
    \end{align*}
    The individual summands consist of nonnegative terms and thus it follows $x_m^*=0$ for $m>2n+1$, as well as
     \begin{align*}
        \lvert x_m^{(k)} -x_m^* \rvert^2< \varepsilon 
    \end{align*}
    for $m=1,\dots,2n+1$.
    It holds $x_m^{(k)} \geq 0$ for $m=1,\dots,2n+1$, and $x_{2i+1}^{(k)} \geq \gamma$ for $i=1,\dots,n$ and for all $k\in \mathbb{N}$. Thus, it follows $x_{2i+1}^* \geq \gamma$ for $i=1,..,n$ and $x_m^{*} \geq 0$ for $m=1,\dots,2n+1$ yielding $x^* \in \mathcal{D}(F_n)$.\\
    Let $\{a_m\},\{b_m\} \in \mathcal{D}(F_n)$ and $\lambda \in [0,1]$. We prove that $$\{c_m\}:= \lambda \{a_m\}+ (1-\lambda)\{b_m\} \in \mathcal{D}(F_n).$$
    Since $a_m=b_m=0$, it follows $c_m=0$ for $m>2n+1$.
    Considering $c_m$ for fixed $m\leq 2n+1$ as a function of $\lambda$, we get $\frac{d c_m}{d \lambda}=a_m-b_m$. For $a_m\geq b_m$ the smallest value for $\lambda=0$ can be found with $c_m=b_m$. For $a_m < b_m$ the minimum is $c_m=a_m$ with $\lambda=1$. In both cases $c_m \geq 0$ holds for all $m \leq 2n+1$ and $c_{2i+1} \geq \gamma$ for all $i=1,\dots,n.$ \
    Thus, $\{c_m\} \in \mathcal{D}(F_n)$ proving the convexity of $\mathcal{D}(F_n)$.
      
\item[(c)] We show that $F_n$ is weak-to-weak sequentially continuous. Let $\{x^{(k)}\}$ be a sequence in $\mathcal{D}(F_n)$ with $x^{(k)}  \rightharpoonup x^* $ in $\l2$ and $x^* \in \mathcal{D}(F)$.
Our aim is to prove the convergence
 $$F_n(x^{(k)}) =:\sigma^{(k)} \rightharpoonup F_n(x^*):=\sigma^*$$ in $L^2([0,T])$.
By the uniform boundedness principle, every weakly convergent sequence is bounded. Moreover, in $\ell^2(\mathbb{N})$, it holds that a sequence is weakly convergent if and only if it is bounded and converges component-by-component.
In $L^2([0,T])$ a sequence converges weakly, if and only if it is bounded and converges pointwise.
We apply this in the following steps:
\begin{enumerate}
\item[1.] We show that $F_n$ maps a bounded sequence 
$\{ x^{(k)} \} \subset \mathcal{D}(F_n)\subset\ell^2(\mathbb{N})$ to a bounded sequence.
\item[2.] We show that from weak convergence in $\ell^2(\mathbb{N})$, and thus boundedness and component-by-component convergence of $x^{(k)}_m \to x^*_m$, the pointwise convergence of $\sigma^{(k)} \xrightarrow{pointwise}   \sigma^*$ in $L^2([0,T])$ follows.
\end{enumerate}
From 1. and 2. we obtain the weak convergence $\sigma^{(k)}\rightharpoonup \sigma^*$.\\

Proof (of 1. and 2.):
\begin{enumerate}[leftmargin=0.5cm]
\item[1.] We have $x^{(k)}=\{x_m^{(k)}\}_{m\in \mathbb{N}} $ for all $k\in \mathbb{N}$. Moreover, $\{ x^{(k)} \}$ is bounded as a weakly convergent sequence, that is, there exists $S\in \mathbb{R}$, such that $\lVert x^{(k)} \rVert_{\ell^2{\mathbb{N}}}= \left( \sum \limits_{m=1}^{\infty} \lvert x_m^{(k)} \rvert^2  \right)^{\frac{1}{2}} \leq S$ holds for all $k\in \mathbb{N}$. Since $x^{(k)} \in \mathcal{D}(F_n)$ we have $x_m^{(k)}=0$ for $m>2n+1$ and $x_m^{(k)}\geq0$ for $m\leq 2n+1$. It follows that $x_m^{(k)}<S$ for all $m,k \in \mathbb{N}$. 
Let $k\in \mathbb{N}$. Then,
\begin{align*}
\lVert  \sigma^{(k)} \rVert_{L^2([0,T])}^2 = &\int \limits_0^{\frac{\bar{\varepsilon}}{\dot{\varepsilon}}} \left[  x_1^{(k)}\dot{\varepsilon} t  + \sum \limits_{j=1}^n  \frac{x^{(k)}_{2j} x^{(k)}_{2j+1} \dot{\varepsilon} }{2} \left(  1-\exp\left(  -\frac{2}{x^{(k)}_{2j+1}}t  \right)  \right) \right]^2 dt\\
 +& \int \limits_{\frac{\bar{\varepsilon}}{\dot{\varepsilon}}}^T  \left[ x^{(k)}_1\bar{\varepsilon} - \sum \limits_{j=1}^n \frac{x^{(k)}_{2j}x^{(k)}_{2j+1}\dot{\varepsilon}}{2} \left( 1-\exp\left( \frac{2\bar{\varepsilon}}{x^{(k)}_{2j+1} \dot{\varepsilon}} \right) \right) \exp\left( -\frac{2}{x^{(k)}_{2j+1}}t \right) \right]^2 dt\\
 &\leq \int \limits_0^T \left[  x^{(k)}_1 \bar{\varepsilon} + \sum \limits_{j=1}^n \frac{x^{(k)}_{2j} x^{(k)}_{2j+1} \dot{\varepsilon}}{2}  \left(   1-\exp\left( -\frac{2\bar{\varepsilon}}{x^{(k)}_{2j+1}\dot{\varepsilon}} \right) \right)  \right]^2 dt \\
 &\leq T \left[  x^{(k)}_1 \bar{\varepsilon} + \sum \limits_{j=1}^n \frac{x^{(k)}_{2j} x^{(k)}_{2j+1} \dot{\varepsilon}}{2}    \right]^2 \\ 
  &\leq T \left[  S \bar{\varepsilon} + \frac{n S^2 \dot{\varepsilon}}{2}    \right]^2 < \infty.
\end{align*}
Thus, $\sigma^{(k)}$ is bounded for all $k\in \mathbb{N}$ implying that the sequence $\{ F_n(x^{(k)}) \}$ is bounded.%\\

\item[2.] We have component-by-component convergence $x_m^{(k)} \xrightarrow{k} x_m^*$ for all $m\in \mathbb{N}$, that means $\forall m\in \mathbb{N} \ \forall \zeta>0 \ \exists K(m) \in \mathbb{N}$, such that $\forall k \geq K(m)$ holds $ | x_m^{(k)}- x_m^*| < \zeta$. Since $x^{(k)} \in \mathcal{D}(F_n), \ k \in \mathbb{N}, \ x^* \in \mathcal{D}(F)$ and thus $x^{(k)}_m=x^*_m=0$ for $m>2n+1, \ k\in \mathbb{N}$, we can choose $K:=\max\{K(m):m=1,..,2n+1\}$ such that $\forall \zeta>0 \ \exists K\in \mathbb{N}$ with $\forall m\in \mathbb{N}, \ \forall k \geq K$ holds $ | x_m^{(k)}- x_m^*| < \zeta$.
We want to show that $\sigma^{(k)} \to   \sigma^*$ pointwise, which means $ \forall t \in [0,T], \ \forall \zeta>0 \ \exists M\in \mathbb{N} $, such that $\forall k \geq M$ holds $ | \sigma^{(k)}(t)- \sigma^*(t) | <\zeta$.
Let $t \in [0,\frac{\bar{\varepsilon}}{\dot{\varepsilon}}]$ and $k\geq K$. We estimate
 \begin{align*}
& |\sigma^{(k)}(t)- \sigma^*(t) |\\
 = &\left| \left( x_1^{(k)} - x_1^* \right) \dot{\varepsilon} t + \frac{\dot{\varepsilon}}{2}\sum \limits_{j=1}^n 
 \left[ x_{2j}^{(k)} x_{2j+1}^{(k)} \left(  1-\exp\left(  -\frac{2}{x_{2j+1}^{(k)}}t  \right)  \right) - 
x_{2j}^* x_{2j+1}^*  \left(  1-\exp\left(  -\frac{2}{x_{2j+1}^*}t  \right)  \right) 
  \right]  \right| \\
 < & \zeta \dot{\varepsilon}t + \frac{\dot{\varepsilon}}{2}\sum \limits_{j=1}^n  \left|  x_{2j}^{(k)} x_{2j+1}^{(k)} \left(  1-\exp\left(  -\frac{2}{x_{2j+1}^{(k)}}t  \right)  \right) - 
x_{2j}^* x_{2j+1}^*  \left(  1-\exp\left(  -\frac{2}{x_{2j+1}^*}t  \right)  \right) 
  \right|  . 
 \end{align*}
For fixed $j\in \{ 1,\dots,n\}$ it follows
\begin{align*}
&\left|  x_{2j}^{(k)} x_{2j+1}^{(k)} \left(  1-\exp\left(  -\frac{2}{x_{2j+1}^{(k)}}t  \right)  \right) - 
x_{2j}^* x_{2j+1}^*  \left(  1-\exp\left(  -\frac{2}{x_{2j+1}^*}t  \right)  \right) \right|\\
  &=\left| x_{2j}^{(k)} x_{2j+1}^{(k)} \left(  1-\exp\left(  -\frac{2}{x_{2j+1}^{(k)}}t  \right)  \right) 
- x_{2j}^{(k)} x_{2j+1}^{(k)} \left(  1-\exp\left(  -\frac{2}{x_{2j+1}^*}t  \right)  \right) \right. \\
&\ \ \ \left. + x_{2j}^{(k)} x_{2j+1}^{(k)} \left(  1-\exp\left(  -\frac{2}{x_{2j+1}^*}t  \right)  \right)
-x_{2j}^* x_{2j+1}^*  \left(  1-\exp\left(  -\frac{2}{x_{2j+1}^*}t  \right)  \right) \right|\\ 
& \leq  \left| x_{2j}^{(k)} x_{2j+1}^{(k)} \left(  \exp\left(  -\frac{2}{x_{2j+1}^*}t  \right)  - \exp\left(  -\frac{2}{x_{2j+1}^{(k)}}t  \right) \right)  \right| \\
& \ \ \ + \left|   \left( x_{2j}^{(k)} x_{2j+1}^{(k)} - x_{2j}^* x_{2j+1}^*  \right) \left(  1-\exp\left(  -\frac{2}{x_{2j+1}^*}t  \right)  \right) \right|\\
&\leq  S^2 \left|   \exp\left(  -\frac{2}{x_{2j+1}^*}t  \right)  - \exp\left(  -\frac{2}{x_{2j+1}^{(k)}}t  \right)   \right|
 + \left|   x_{2j}^{(k)} x_{2j+1}^{(k)} - x_{2j}^* x_{2j+1}^*   \right|.
\end{align*}
We consider the two parts separately, starting with the second.
\begin{align*}
\left|    x_{2j}^{(k)} x_{2j+1}^{(k)} - x_{2j}^* x_{2j+1}^* \right|
&=\left|    x_{2j}^{(k)} x_{2j+1}^{(k)} -  x_{2j}^{(k)} x_{2j+1}^*  + x_{2j}^{(k)} x_{2j+1}^* - x_{2j}^* x_{2j+1}^* \right| \\
&\leq   x_{2j}^{(k)}  \left|   x_{2j+1}^{(k)} -  x_{2j+1}^*  \right| + x_{2j+1}^* \left| x_{2j}^{(k)}  - x_{2j}^*  \right|
\leq 2S \zeta.
\end{align*}
To prove the estimate for $\left| \exp\left( -\frac{2t}{x_{2j+1}^{(k)}}  \right) - \exp\left( -\frac{2t}{x_{2j+1}^*} \right) \right| $, we exploit that the exponential function is uniformly continuous on the compact interval $J_1:=[-\frac{2T}{\gamma},0]$. So, for any given $\eta > 0$ we find some $\delta=\delta(\eta)$ such that for any $x^{(1)}, x^{(2)}\in J_1$ with $|x^{(1)}-x^{(2)}|< \delta$ we have $\left|\exp(x^{(1)})-\exp(x^{(2)})\right|<\eta$.
Computing 
\begin{align*}
|x^{(1)}-x^{(2)}|
= \left|   \frac{2t}{x_{2j+1}^*}    -\frac{2t}{x_{2j+1}^{(k)}}  \right|
= \frac{2t\left| x_{2j+1}^{(k)} - x_{2j+1}^*  \right|}{x_{2j+1}^{(k)}x_{2j+1}^{*}} < \frac{2T}{\gamma^2} \zeta \end{align*}
and choosing $\zeta$ sufficiently small, such that $(2T\zeta/\gamma^2)<\delta (\eta)$, we conclude 
$$\left| \exp\left( -\frac{2t}{x_{2j+1}^{(k)}}  \right) - \exp\left( -\frac{2t}{x_{2j+1}^*} \right) \right| <\eta$$ 
for an arbitrarily given $\eta>0$. This leads to 
\begin{align*}
    |\sigma^{(k)}(t)- \sigma^*(t) | \leq  \zeta \dot{\varepsilon}t + \frac{\dot{\varepsilon}}{2}\sum \limits_{j=1}^n \left(S^2 \eta + 2S \zeta \right)
\end{align*}
and we obtain pointwise convergence $\sigma^{(k)} \to  \sigma^*$ for $t\in [0,\frac{\bar{\varepsilon}}{\dot{\varepsilon}}]$.
Next we show pointwise convergence in $[\frac{\bar{\varepsilon}}{\dot{\varepsilon}},T]$.
To this end, let $k\geq K$, $t \in [\frac{\bar{\varepsilon}}{\dot{\varepsilon}},T]$ and 
\begin{align*}
& |\sigma^{(k)}(t)- \sigma^*(t) |\\
 = &\left| \left( x_1^{(k)} - x_1^* \right) \bar{\varepsilon} + \frac{\dot{\varepsilon}}{2}\sum \limits_{j=1}^n 
 \left[ x_{2j}^{(k)} x_{2j+1}^{(k)} \left(  1-\exp\left( \frac{2\bar{\varepsilon}}{x_{2j+1}^{(k)}\dot{\varepsilon}}  \right)  \right)\exp\left(  -\frac{2t}{x_{2j+1}^{(k)}}  \right)  \right. \right. \\
 &\hspace{3.8cm} \left. \left.
-x_{2j}^* x_{2j+1}^*  \left(  1-\exp\left( \frac{2\bar{\varepsilon}}{x_{2j+1}^*\dot{\varepsilon}}  \right)  \right)\exp\left(  -\frac{2t}{x_{2j+1}^*}  \right)
  \right]  \right|.
 \end{align*} 
This yields $\left| \left( x_1^{(k)} - x_1^* \right) \right| \leq \zeta \bar{\varepsilon}$ and, for fixed $j\in \{ 1,\dots,n\}$, we obtain the estimate for the individual summands:
\begin{align*}
&\left| \left[ x_{2j}^{(k)} x_{2j+1}^{(k)} \left(  1-\exp\left( \frac{2\bar{\varepsilon}}{x_{2j+1}^{(k)}\dot{\varepsilon}}  \right)  \right)\exp\left(  -\frac{2t}{x_{2j+1}^{(k)}}  \right) 
-x_{2j}^* x_{2j+1}^*  \left(  1-\exp\left( \frac{2\bar{\varepsilon}}{x_{2j+1}^*\dot{\varepsilon}}  \right)  \right)\exp\left(  -\frac{2t}{x_{2j+1}^*}  \right)
 \right]  \right| \\
&=\left| x_{2j}^{(k)} x_{2j+1}^{(k)} \left(  1-\exp\left( \frac{2\bar{\varepsilon}}{x_{2j+1}^{(k)}\dot{\varepsilon}}  \right)  \right)\exp\left(  -\frac{2t}{x_{2j+1}^{(k)}}  \right) 
- x_{2j}^{(k)} x_{2j+1}^{(k)} \left(  1-\exp\left( \frac{2\bar{\varepsilon}}{x_{2j+1}^{(k)}\dot{\varepsilon}}  \right)  \right)\exp\left(  -\frac{2t}{x_{2j+1}^*}  \right)  \right. \\
&\ \ \  \left. +x_{2j}^{(k)} x_{2j+1}^{(k)} \left(  1-\exp\left( \frac{2\bar{\varepsilon}}{x_{2j+1}^{(k)}\dot{\varepsilon}}  \right)  \right)\exp\left(  -\frac{2t}{x_{2j+1}^*}  \right) 
- x_{2j}^* x_{2j+1}^*  \left(  1-\exp\left( \frac{2\bar{\varepsilon}}{x_{2j+1}^*\dot{\varepsilon}}  \right)  \right)\exp\left(  -\frac{2t}{x_{2j+1}^*}  \right) \right|\\
&\leq \left| x_{2j}^{(k)} x_{2j+1}^{(k)}  \left(  1-\exp\left( \frac{2\bar{\varepsilon}}{x_{2j+1}^{(k)}\dot{\varepsilon}}  \right)  \right) \left(  \exp\left(  -\frac{2t}{x_{2j+1}^{(k)}}  \right)  - \exp\left(  -\frac{2t}{x_{2j+1}^*}  \right)  \right)  \right| \\
&\ \ \  + \left|   \left[   x_{2j}^{(k)} x_{2j+1}^{(k)} \left(  1-\exp\left( \frac{2\bar{\varepsilon}}{x_{2j+1}^{(k)}\dot{\varepsilon}}  \right)  \right) -   x_{2j}^* x_{2j+1}^*  \left(  1-\exp\left( \frac{2\bar{\varepsilon}}{x_{2j+1}^*\dot{\varepsilon}}  \right)  \right)  \right] \exp\left(  -\frac{2t}{x_{2j+1}^*}  \right)  \right| \\
&\leq S^2 \left|  \exp\left(  -\frac{2t}{x_{2j+1}^{(k)}}  \right)  - \exp\left(  -\frac{2t}{x_{2j+1}^*}  \right)    \right| 
 + \left|  \left( x_{2j}^{(k)} x_{2j+1}^{(k)} - x_{2j}^* x_{2j+1}^*  \right) \left(  1-\exp\left( \frac{2\bar{\varepsilon}}{x_{2j+1}^{(k)}\dot{\varepsilon}}  \right)  \right) \right|  \\
&\ \ \  + \left|   x_{2j}^* x_{2j+1}^* \left( \exp\left( \frac{2\bar{\varepsilon}}{x_{2j+1}^*\dot{\varepsilon}}  \right)- \exp\left( \frac{2\bar{\varepsilon}}{x_{2j+1}^{(k)}\dot{\varepsilon}}  \right) \right) \right|.
\end{align*}  
We know that 
$$\left| \exp\left( -\frac{2t}{x_{2j+1}^{(k)}}  \right) - \exp\left( -\frac{2t}{x_{2j+1}^*} \right) \right|<\eta$$ 
holds for $J_1:=[-\frac{2T}{\gamma},0]$. 
Concerning the estimate of the second summand we conclude as in the previous proof 
\begin{align*}
& \left|  \left( x_{2j}^{(k)} x_{2j+1}^{(k)} - x_{2j}^* x_{2j+1}^*  \right) \left(  1-\exp\left( \frac{2\bar{\varepsilon}}{x_{2j+1}^{(k)}\dot{\varepsilon}}  \right)  \right) \right|    \\
 &\  \leq 2S  \left| 1-\exp\left( \frac{2\bar{\varepsilon}}{x_{2j+1}^{(k)}\dot{\varepsilon}}  \right) \right| \zeta
\leq 2S \left| 1-\exp\left( \frac{2\bar{\varepsilon}}{\gamma\dot{\varepsilon}}  \right) \right| \zeta \\
   &\  \leq 2S  \exp\left( \frac{2\bar{\varepsilon}}{\gamma\dot{\varepsilon}}  \right) \zeta.
\end{align*}
Again, we exploit the uniform continuity of the exponential function on $J_2:=[-\frac{2\bar{\varepsilon}}{\gamma \dot{\varepsilon}},0]$ to obtain $$|x^{(1)}-x^{(2)}|
= \left| \frac{2\bar{\varepsilon}}{x_{2j+1}^* \dot{\varepsilon}}    -\frac{2\bar{\varepsilon}}{x_{2j+1}^{(k)}\dot{\varepsilon}}  \right|
< \frac{2\bar{\varepsilon}}{\gamma^2\dot{\varepsilon}} \zeta,$$ 
from what
$$\left| \exp\left( \frac{2\bar{\varepsilon}}{x_{2j+1}^*\dot{\varepsilon}} \right)- \exp\left( \frac{2\bar{\varepsilon}}{x_{2j+1}^{(k)}\dot{\varepsilon}}  \right)\right| <\eta$$ 
for an arbitrarily given $\eta>0$ follows, if only $\zeta$ is sufficiently small. Finally,
\begin{align*}
    |\sigma^{(k)}(t)- \sigma^*(t) | \leq  \zeta \bar{\varepsilon}+ \frac{\dot{\varepsilon}}{2}\sum \limits_{j=1}^n \left (S^2 \eta + 2S \exp\left(\frac{2 \bar{\varepsilon}}{\gamma \dot{\varepsilon}} \right)\zeta + S^2 \eta \right)
\end{align*}
holds true proving the pointwise convergence $\sigma^{(k)} \to   \sigma^*$ pointwise for $t\in [\frac{\bar{\varepsilon}}{\dot{\varepsilon}},T]$ and, thus, for all $t\in [0,T]$. \\
\qed
\end{enumerate}
\end{enumerate}

These investigations imply the first block of assumptions \ref{Ann: 1st block}. 

Next, we turn to the second block of assumptions \ref{Ann: 2nd block} (applied to \eqref{eq: Tikh in Bayes algo section}). The exponent in \eqref{eq: Tikh in Bayes algo section} is $p=2$ and, thus, \ref{Ann: 2nd block}(a) is satisfied.
Part (b) and (c) depend on the used penalty term $\Omega(x)$. For the numerical evaluations in section \ref{sec: regularization} we apply the penalties
\begin{align*}
   \Omega_1(x)=\frac{1}{2}\lVert x\rVert ^2 \qquad\text{and}\qquad \Omega_2(x)= \frac{1}{2} x_{3}^2
\end{align*}
for $x \in \l2$ and $0< \gamma \leq x_3 \in \mathbb{R}$. 
Both penalties, $\Omega_1$ as well as $\Omega_2$ have a nonempty domain and $ \mathcal{D}:= \mathcal{D}(F) \cap \mathcal{D}(\Omega_i) \neq \emptyset$ holds true for $i=1,2$. Both functionals are convex and continuous, which also implies the lower semi-continuity. Thus, assumption \ref{Ann: 2nd block}(b) is satisfied for $\Omega_1$ and $\Omega_2$. As already mentioned, penalty terms $\Omega (x)$ of the form \eqref{eq: omega norm form} satisfy assumption \ref{Ann: 2nd block} (c). So, (c) is fulfilled for $\Omega_1$, $\Omega_2$.
The existence of a $\Omega$-minimizing solution $x^+$ to the problem \eqref{eq:F_n equation} which belongs to the Bregman domain $\mathcal{D}_B(\Omega)$ is postulated such that \ref{Ann: 2nd block}(d) is also satisfied.

In this sense theorems \ref{theorem:existence} - \ref{theorem: convergence} can be applied to \eqref{eq:F_n equation} for fixed $n \in I$ and the associated Tikhonov functional \eqref{eq: Tikh in Bayes algo section}.

%%%%%%%%%%%%%%%%%%%%%%%%%%%%%%%%%%%%%%%%%%%%%%%%%%%%%%%%%%%%%%%%

%\textbf{Application to the overall method}\\
In the following we aim to prove the regularization property of algorithm \ref{alg:Bayes}. The next defintion is adapted from the well-known definition of a $\Omega$-minimizing solution (c.f. definition \ref{dfn: Omega min solution}).

\begin{dfn}
We call $(n^+,x^+)\in \mathcal{D}(F) \subseteq \mathbb{N}\times\l2 $ a $\logg$-$\Omega$-\textit{minimizing solution} of \eqref{eq:unsere_operatorgleichung} if 
\begin{align*}
    F(n^+,x^+)=\sigma 
\end{align*}
and
\begin{align*}
    \logg(n^+)+\Omega(x^+)=\inf \{ \logg(\tilde{n}) + \Omega(\tilde{x}) \, : \, (\tilde{n},\tilde{x})\in \mathcal{D}(F), \, F(\tilde{n},\tilde{x})=\sigma \}
\end{align*}
holds true for exact data $\sigma \in L^2 ([0,T])$ in \eqref{eq:unsere_operatorgleichung}.
\end{dfn}

Since $\mathbb{N}\times \ell^2 (\mathbb{N})$ is a carteian product of a semigroup and a Hilbert space we extend the definiton of a regularization method to such structures.

\begin{dfn}[regularization in $\mathbb{N}\times \ell^2 (\mathbb{N})$]
\label{def:regularization_semigroups}
A mapping that assigns each pair $(\sigma^\delta,\alpha) \in L^2([0,T]) \times (0,\bar{\alpha})$ with $0< \bar{\alpha} \leq + \infty$ to well-defined elements $(n^\delta_\alpha,x^\delta_\alpha) \in \mathcal{D}(F)$ is called \textit{regularization (regularization method)} for \eqref{eq:unsere_operatorgleichung} if there is a suitable choice $\alpha=\alpha(y^\delta,\delta)$ of the regularization parameter , such that for any sequence $\{ \sigma_m \}_{m=1}^{\infty}\subset L^2([0,T])$ with $$\lVert \sigma_m-\sigma \rVert \leq \delta_m$$ and $\delta_m \to 0$ for $m \to \infty$, the corresponding regularized solutions $(n_{\alpha(y_m,\delta_m)}^{\delta_m},x_{\alpha(y_m,\delta_m)}^{\delta_m})$ converge in a well-defined sense to a solution $(n^+,x^+)$ of \eqref{eq:unsere_operatorgleichung}. Since the latter is in general not unique, regularized solutions must converge to $\logg$-$\Omega$-minimizing solutions of \eqref{eq:unsere_operatorgleichung}. In the case of nonuniqueness, different subsequences of regularized solutions may converge to different $\logg$-$\Omega$-minimizing solutions of \eqref{eq:unsere_operatorgleichung}.\\
\end{dfn}

For convergence, we consider primarily the product topologies $(\tau_d, \tau_0)$ and $(\tau_d, \tau_w)$.

%%%%%%%%%%%%%%%%%%%%%%%
%\paragraph{Existence of a minimizer of \eqref{eq: T_alpha(n,x)}}
%%%%%%%%%%%%%%%%%%%%%%%

We have shown that for all $\alpha>0$, $n \in I$ and $\sigma^\delta \in L^2([0,T])$, there exists a regularized solution $x_\alpha^\delta\in \mathcal{D}(F_n)$ minimizing the functional $T_{\alpha,n}(x)$ \eqref{eq: Tikh in Bayes algo section} in $\mathcal{D}(F_n)$.
Since $I$ is finite, then there exists a solution $(n_\alpha^\delta,x_\alpha^\delta)$ that minimizes the functional $T_\alpha(n,x)$\eqref{eq:Tikhonov_Bayes} in $(n,x)\in \mathcal{D}(F)$.

%%%%%%%%%%%%%%%%%%%%%%%%
%\paragraph{Stability}
%%%%%%%%%%%%%%%%%%%%%%%%

\begin{theorem}[Stability]
\label{theorem: stability overall}
For all $\alpha>0$, the minimizers of \eqref{eq:Tikhonov_Bayes} are stable with respect to the data $\sigma^\delta$. Let $\{ \sigma_m \}$ be a sequence of data with $\lim \limits_{n\to\infty}\lVert \sigma_m-\sigma^\delta \rVert_{L^2([0,T])}=0$ and $\{ (n_m, x_m)\}$ be the corresponding sequence of minimizers of 
\begin{align*}
    \min_{(n,x) \in \mathcal{D}(F)} \left \{ \frac{1}{2}\lVert F(n,x) - \sigma_m \rVert^2 +\alpha \left[ \logg(n)+\Omega(x)\right] \right \}.
\end{align*}
Then, $\{( n_m, x_m )\}$ has a $(\tau_d,\tau_w)$-convergent subsequence $\{( n_{m_k}, x_{m_k}) \}$.  The limit of each such subsequence is a minimizer $(n_\alpha^\delta, x_\alpha^\delta)$ of \eqref{eq:Tikhonov_Bayes}.
Moreover, for any such $(\tau_d,\tau_w)$-convergent subsequence 
$$\lim \limits_{k\to \infty} \left( \logg(n_{m_k}) + \Omega(x_{m_k}) \right)=\logg(n_\alpha^\delta) + \Omega(x_\alpha^\delta)$$
holds true. 
\end{theorem}

\textit{Proof:}\\
Since $\{n_m\} \subset I $ and finite $I$, there must be a $\tau_d$-convergent subsequence $\{ n_{m_k} \}$, that means, there exists a $k_0 \in \mathbb{N}$ such that for all $k \geq k_0$ we have $n_{m_k}=\tilde{n}$.
For $k \geq k_0$, $x_{m_k}$ can be considered as the minimizer of 
\begin{align*}
\frac{1}{2}\lVert F_{\tilde{n}}(x) - \sigma_m \rVert^2 +\alpha \Omega(x)
\end{align*}
for fixed $\tilde{n}$. 
Then, by theorem \ref{theorem: stability}, there exists a subsequence $\{x_{m_k}\} $ with $x_{m_k} \rightharpoonup \tilde{x}$, where $\tilde{x}$ is a minimizer of 
$$T_{\alpha,\tilde{n}}(x)= \frac{1}{2}\lVert F_{\tilde{n}}(x ) - \sigma^\delta \rVert^2 +\alpha \Omega(x)$$ 
and 
$$\lim \limits_{k\to \infty} \Omega(x_{m_k})=\Omega(x_\alpha^\delta).$$
Moreover, $F(n_{m_k},x_{m_k})= F_{\tilde{n}}(x_{m_k}) \rightharpoonup F_{\tilde{n}}(\tilde{x}) $ as $k \to \infty$, since $F_{\tilde{n}}$ is weak-to-weak sequentially continuous and hence it follows that $F(n_{m_k},x_{m_k}) - \sigma_{m_k} \rightharpoonup F(\tilde{n}, \tilde{x}) - \sigma^\delta $ for $k\to \infty$. 
Since the $L^2$-norm and $\Omega$ are weakly lower semi-continuous, we deduce 
\begin{align*}
    \frac{1}{2} \lVert F(\tilde{n},\tilde{x})- \sigma^\delta \rVert^2
    &\leq \liminf \limits_{k\to \infty} \frac{1}{2} \lVert F(n_{m_k},x_{m_k}) - \sigma_{m_k} \rVert ^2, \\
    \Omega(\tilde{x}) &\leq \liminf \limits_{k \to \infty} \Omega(x_{m_k}). 
\end{align*}
Since for $k \geq k_0$ we have $\logg (n_{m_k}) =\logg(\tilde{n})$, it furthermore follows $\logg(\tilde{n}) \leq \liminf \limits_{k \to \infty}  \logg(n_{m_k}) $ and thus \begin{align*}
  \hspace{-1cm} \frac{1}{2} \lVert F(\tilde{n},\tilde{x})- \sigma^\delta \rVert^2 + \alpha\left[\logg(\tilde{n})+ \Omega(\tilde{x})   \right]
    &\leq \liminf \limits_{k\to \infty} \left( \frac{1}{2} \lVert F(n_{m_k},x_{m_k}) - \sigma_{m_k} \rVert ^2 + \alpha \left[ \logg (n_{m_k}) +  \Omega(x_{m_k})   \right] \right)\\
    &\leq
    \lim \limits_{k \to \infty} \frac{1}{2}\lVert F(n,x) - \sigma_{m_k} \rVert^2 +\alpha \left[\logg (n) + \Omega(x) \right] \\
    &= \frac{1}{2}\lVert F(n,x) - \sigma^\delta \rVert^2 + \alpha \left[\logg (n) + \Omega(x) \right],
\end{align*}
where $(n,x)\in \mathcal{D}(F)$ is arbitrary. This shows that $(\tilde{n},\tilde{x})$ is a minimizer of \eqref{eq:Tikhonov_Bayes}. 
If we choose $x=\tilde{x}$ and $n= \tilde{n}$ on the right-hand side, it follows \begin{align*}
   \hspace{-0.5cm}   \frac{1}{2} \lVert F(\tilde{n},\tilde{x})- \sigma^\delta \rVert^2 + \alpha\left[\logg(\tilde{n})+ \Omega(\tilde{x})   \right]
    = \lim \limits_{k\to \infty} \left( \frac{1}{2} \lVert F(n_{m_k},x_{m_k}) - \sigma_{m_k} \rVert ^2 + \alpha \left[ \logg (n_{m_k}) +  \Omega(x_{m_k})   \right] \right). 
\end{align*}
Moreover, it holds $\lim \limits_{k\to \infty} \Omega(x_{m_k})=\Omega(x_\alpha^\delta)$ and $\lim \limits_{k\to \infty} \logg(n_{m_k})=\logg(n_\alpha^\delta)$. Thus, we finally get 
$$\lim \limits_{k\to \infty} \left( \logg(n_{m_k}) + \Omega(x_{m_k}) \right)=\logg(n_\alpha^\delta) + \Omega(x_\alpha^\delta).$$ 
\qed\\

%\paragraph{Convergence}

We finally prove convergence of the method and, this way, the regularization property.

\begin{theorem}[Convergence]
\label{theorem: convergence overall}
Let $\{ \sigma_m:= \sigma^{\delta_m} \} \subset L^2([0,T])$ be a sequence of perturbed data, $\sigma\in F(\mathcal{D}(F))$ exact data with $\lVert \sigma - \sigma_m \rVert \leq \delta_m $ for a sequence $\{ \delta_m >0 \}$ of noise levels converging monotonically to zero. Moreover, we consider a sequence $\{ \alpha_m >0 \}$ of regularization parameters and an associated sequence $\{( n_m:=n_{\alpha_m}^{\delta_m} , x_m:=x_{\alpha_m}^{\delta_m}) \}$ of regularized solutions which are minimizers of 
\begin{align*}
    \frac{1}{2}\lVert F(n,x) - \sigma_m \rVert^2 +\alpha_m \left[ \logg(n) + \Omega(x) \right] \text{ for } (n,x) \in \mathcal{D}(F).
\end{align*}
Under the conditions 
\begin{align}
\label{eq: convergence omegaloggBed}
    &\limsup \limits_{m\to \infty} \left(\logg (n_m) + \Omega(x_m) \right) \leq \logg(n_0)+ \Omega(x_0) \nonumber \\
    \text{ for all } &(n_0,x_0) \in \{ n\in I, x\in \mathcal{D}:= \mathcal{D}(F_{n})\cap \mathcal{D}(\Omega) : F(n,x)=\sigma \},
\end{align}
and 
\begin{align}
\label{eq: convergence condition 2}
    \lim \limits_{m\to \infty} \lVert F(n_m,x_m)-\sigma_m \rVert =0 
\end{align}
the sequence $\{( n_m,x_m) \}$ has a $(\tau_d,\tau_w)$-convergent subsequence, where each limit is a $\logg$-$\Omega$-minimizing solution $(n^+,x^+) \in \mathcal{D}(F)$ of \eqref{eq:unsere_operatorgleichung}. Additionally, if the $\logg$-$\Omega$-minimizing solution $(n^+,x^+)$ is unique, we obtain $(\tau_d,\tau_w)$-convergence $$(n_m,x_m) \to (n^+,x^+) \text{ in } \mathbb{N}\times \l2.$$ 
\end{theorem}

\textit{Proof:}\\
As in the proof of theorem \ref{theorem: stability overall}, we can conclude for the sequence $\{ n_m \} $ that it has a convergent subsequence $\{n_{m_k} \}$ with $n_{m_k} \to \tilde{n}$. That is, there exists a $k_0\in \mathbb{N}$ such that for all $k\geq k_0$ we have $n_{m_k}=\tilde{n}$.
From $\eqref{eq: convergence condition 2}$ it follows from $\delta_m \to 0$ for $m \to +\infty$ that
\begin{align}
\label{eq: convergence proof triangle ungl}
    \lim \limits_{m\to \infty} \lVert F(n_m,x_m)- \sigma \rVert 
    \leq \lim \limits_{m\to \infty}  \lVert F(n_m,x_m)- \sigma_m \rVert+ \lVert \sigma_m - \sigma \rVert =0 .
\end{align}
Additionally 
\begin{align*}
    \lVert F(n_m,x_m) \rVert \leq \lVert \sigma_m \rVert + \lVert F(n_m,x_m)-\sigma_m \rVert 
\end{align*}
holds true for all $m \in \mathbb{N}$. Since $\lim \limits_{m\to \infty}\sigma_m=\sigma$ and $\lim \limits_{m\to \infty} \lVert F(n_m,x_m)-\sigma_m \rVert=0 $ according to \eqref{eq: convergence condition 2}, both $\lVert \sigma_m \rVert$ and $\lVert F(n_m,x_m)-\sigma_m \rVert$ are bounded. Thus $\{ F(n_m,x_m) \}$ is bounded and there exists a weakly convergent subsequence $\{ F(n_{m_k},x_{m_k}) \}$ in $L^2([0,T])$. So we can choose a subsequence $\{ (n_{m_k},x_{m_k} )\}$ such that $n_{m_k}=\tilde{n}$ is constant for all $k\in \mathbb{N}$ yielding 
$$\{F(n_{m_k},x_{m_k}) \}=\{F(\tilde{n},x_{m_k}) \}. $$ 
By theorem \ref{theorem: convergence} there exists a subsequence $\{x_{m_k}\}$ with $x_{m_k} \rightharpoonup \tilde{x}$ and \eqref{eq: convergence proof triangle ungl} leads to 
$$F(\tilde{n},\tilde{x})=\sigma.$$
From condition \eqref{eq: convergence omegaloggBed} we deduce
\begin{align*}
    \logg(\tilde{n}) + \Omega(\tilde{x}) \leq \liminf \limits_{k\to \infty}  \left( \logg(n_{m_k}) + \Omega(x_{m_k}) \right) \leq \limsup \limits_{k\to \infty} \left( \logg(n_{m_k}) + \Omega(x_{m_k}) \right)&\\
    \leq \logg(n^+) + \Omega(x^+) \leq \logg(n_0) + \Omega(x_0)& 
\end{align*}
for all $(n_0,x_0)\in \mathcal{D}(F)$ with $F(n_0,x_0)=\sigma$.  
If we set $(n_0,x_0)=(\tilde{n},\tilde{x})$ this leads to $$\logg(\tilde{n})+ \Omega(\tilde{x})=\logg(n^+)+\Omega(x^+)$$ and thus $(\tilde{n},\tilde{x})$ is a $\logg$-$\Omega$-minimizing solution. Moreover, 
$$\lim \limits_{m\to \infty} \left( \logg (n_m)+ \Omega(x_m) \right)= \logg(n^+)+\Omega(x^+).$$ 
If $(n^+,x^+)$ is unique, then we have $(\tau_d,\tau_w)$-convergence to $(n^+,x^+)$ for every subsequence of $(n_m,x_m)$ and, thus, $(\tau_d,\tau_w)$-convergence  
$$(n_m,x_m) \rightharpoonup (n^+,x^+).$$
\qed

%%%%%%%%%%%%%%%%%%%%%%%%%%%%%%%%%%%%%%%%%%%%%%%%%%%%%%

\section{Numerical validation} \label{chap:Numeric}
%mit Vergleich zum Cluster-Algorithmus)
%%%%%%%%%%%%%%%%%%%%%%%%%%%%%%%%%%%%%%%%%%%%%%%%%%%%%%

In this section, we validate the algorithm \ref{alg:Bayes} by checking its performance and comparing it with the cluster algorithm published in \cite{Rothermel2022} and a simple least squares mthod. At first we briefly recall and sketch the cluster algorithm.

% % % % % % % % % % % % % % % % % % % % % % % %
\subsection{The cluster algorithm} \label{sec: cluster}
%%%%%%%%%%%%%%%%%%%%%%%%%%%%%%%%%%%%%%%%%%%%%%%

Note, that the forward operator of the underlying inverse problem is not injective, that is, there exists parameters $(n_1,x_1)$ and $(n_2,x_2)$ with $(n_1,x_1)\neq (n_2,x_2)$ and
$$F\left(n_1,x_1 \right)=F(n_2,x_2).$$

We overcome this non-uniqueness by the following requirement. Assume that the relaxation times are located in different decades (cf. \cite{Goldschmidt2015ModelingBonds}). For example, if $\tau_l \in [10,100]$ for some $l \in \{ 1,\dots , n \}$, then $\tau_j \notin [10,100]$ must hold true for all $j\in \{1, \dots,l-1,l+1, \dots n \}$.
This information is essential for the development of the cluster algorithm.

Let us fix a maximum number of Maxwell elements $N \in \mathbb{N}$ such that the unknown number of Maxwell elements $n^*$ satisfies $n^*\leq N$. Thus, the set of possible numbers of Maxwell elements is represented by $I:=\lbrace 1, \dots, N\rbrace$. The material parameters 
$ x=\big(\mu, \mu_1, \tau_1,\ldots,\mu_N, \tau_N, 0,0,\ldots \big)$ are then computed by a minimization algorithm that works in the same way as the minimization process $R:I \subset \mathbb{N} \to \mathcal{D}(F_n) \subset \l2$, $n \mapsto x$ described in section \ref{sec: Bayes Algo} with the difference that the number of Maxwell elements is now known and given by $N$. This way we obtain a minimizer $x^*$. 

After the minimization process we apply an algorithm to 
$$x^*=\big(\mu^*, \mu_1^*, \tau_1^*,\ldots,\mu_N^*, \tau_N^*, 0,0,\ldots \big)$$ 
that clusters the relaxation times $(\tau_1^*,\ldots,\tau_N^*)$ according to the decade condition. The cluster algorithm consists of two parts:
\begin{enumerate}
\item minimization with $N$ Maxwell elements, output: $x^*$
\item clustering of relaxation times $(\tau_1^*,\ldots,\tau_N^*)$ which reduces $N$ Maxwell elements to $n^*\leq N$ elements.
\end{enumerate}

Let us consider the clustering step 2.) in more detail.
For a decade $[10^k, 10^{k+1}]$ with some $k\in \mathbb{N}_0$, we arrange the Maxwell elements by defining index sets
\begin{align*}
J_k:=\big\{ j \in \lbrace 1,\dots, N \rbrace: \tau_j \in [10^k, 10^{k+1}] \big\},
\end{align*}
where $k\in \mathbb{N}_0$ denotes the decade in which the relaxation time is classified. For the choice of decades it is useful to take into account the physical conditions of the experiment and the material under consideration. 
The maximum number of Maxwell elements $N$ should not be larger than the number of available decades. Accordingly, it makes sense to set $N$ equal to the number of decades.

As next step, the pairs $(\mu_j,\tau_j)$, $j\in J_k$, must be assigned to a Maxwell element. Collecting the relaxation times and thus the Maxwell elements in index sets, the number of non-empty index sets yields the current number $n$ of elements in the material. To cut down multiple Maxwell elements to a single one, we use the following approximate calculation for the new material parameters $(\tilde{\mu}_k, \tilde{\tau}_k)$ with $\tilde{\tau}_k \in [10^k, 10^{k+1}]$ from the already reconstructed material parameters:
 \begin{align}
    \tilde{\mu}_k:= \sum_{j\in J_k} \mu_j^*, \qquad \tilde{\tau}_k:= \sum_{j\in J_k} \frac{\mu_j^*}{\tilde{\mu}_k} \tau_j^*.
    \label{eq:cluster calculation}
    \end{align}
%The disadvantage of this calculation is that there is no analytical way to combine the parameters in such a way that exactly the same stress curve is obtained when the forward operator is applied.
This way, we get the following algorithm: 

\begin{alg}\ \\
\label{alg: cluster}
Input: $N$, $ \sigma \in L^2([0,T])$.
\begin{itemize}
    \item[1)] Compute $x^*$ by the minimization procedure $R(N)= x^*$. 
    \item[2)] Determine the index sets $J_k \neq \emptyset$ for $k\in \mathbb{N}_0$ with $
J_k:=\lbrace j \in \lbrace 1,\dots, N \rbrace: \tau_j^* \in [10^k, 10^{k+1}] \rbrace$ and set $\tilde{n}$ as the number of nonempty index sets.
     \item[3)] For each $k\in \mathbb{N}_0$ with $J_k \neq \emptyset$ compute
     $$\tilde{\mu}_k:= \sum_{j\in J_k} \mu_j^*, \qquad \tilde{\tau}_k:= \sum_{j\in J_k} \frac{\mu_j^*}{\tilde{\mu}_k} \tau_j^*$$
     and set 
     $\tilde{x}:=\left(\mu^*, \tilde{\mu}_1, \tilde{\tau}_1,\dots,\tilde{\mu}_{\tilde{n}}, \tilde{\tau}_{\tilde{n}}, 0,0,\dots \right)$.
     
\end{itemize}
Output: $(\tilde{n},\tilde{x})$.
\end{alg}

In the cluster algorithm, the set of possible numbers of Maxwell elements is  fixed by the maximum number $N$ and not a free parameter to be computed additionally and allowing a specific choice adapted to the given experiment. 
Furthermore, the cluster algorithm requires the classification of relaxation times in decades. This assumption is not necessary for the Bayes algorithm. On the other hand, the Bayes algorithm allows a simple control of the preferred number of Maxwell elements via the success probability $q$ in the prior. Thus, the Bayes algorithm has a more significant regularizing effect, since its penalty term is physically meaningful and derived by statistical inversion theory. 
This will be evident in the following reconstruction results that use perturbed data.

%%%%%%%%%%%%%%%%%%%%%%%%%%%%%%%%%%%
\subsection{Reconstructions from exact data}
%%%%%%%%%%%%%%%%%%%%%%%%%%%%%%%%%%%

\begin{figure}
     \centering
    \subfloat{\def\svgwidth{0.5\textwidth}
     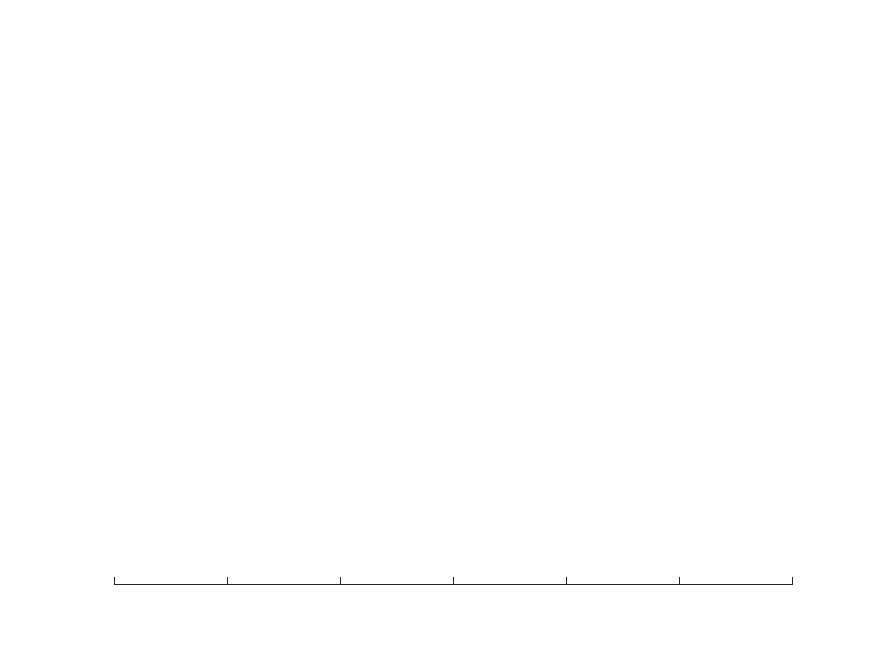}
 \subfloat{\def\svgwidth{0.5\textwidth}
   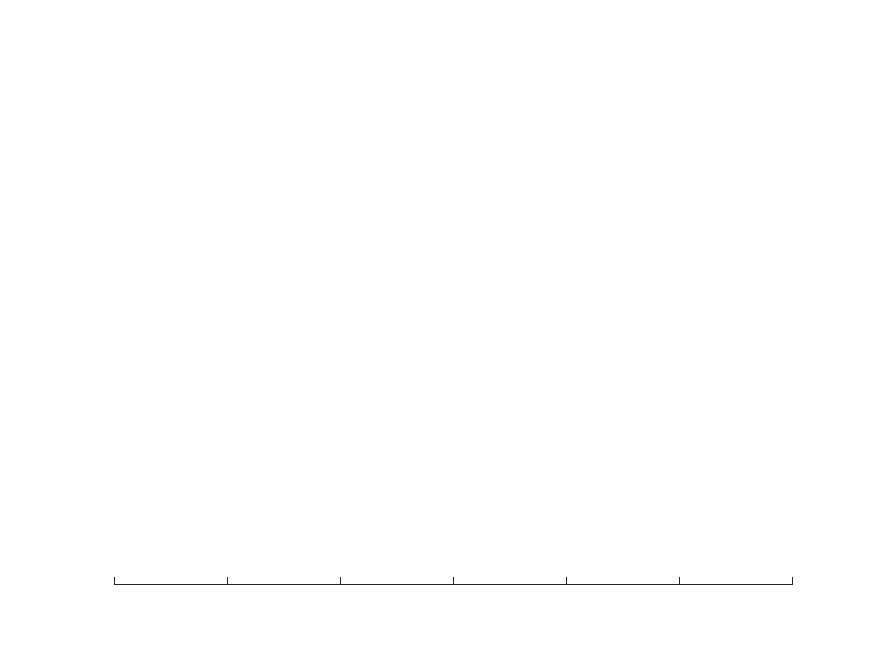}
   \caption{Strain- and stress-data for different displacement rates $\dot{\varepsilon}_u$}
        \label{abb: strain and stress}
\end{figure}

We generate simulated data for given material parameters that serve as the basis of our experiments. In figure \ref{abb: strain and stress} strains for two different displacement rates $\dot{\varepsilon}_u$ and the maximum strain $\bar{\varepsilon}=20\,\%$ are depicted. The fast displacement rate is $10$ mm/s, with the maximum strain reached after $2$ seconds. At the slow displacement rate of $1$ mm/s, the strain of $20 \%$ is achieved only after $20$ seconds. Here, we set $T=100$ seconds.
For the first experiment we choose a material which is characterized by  $n^*=3$ Maxwell elements and parameters $x^*$ as given in table \ref{table:Exact parameters}. 

\begin{table}[h]
\centering
		\begin{tabular}{lcccc}\toprule
			$j$&$0$ &$1$&$2$&$3$ \\ \midrule
			 $\mu^*_j$ [MPa] &10&4&7&1 \\ \midrule 
			 $\tau^*_j$ [s] &-&0.2&3.7&25 \\ \bottomrule
		\end{tabular}
	\caption{Exact material parameters $x^*$ in experiment 1 to simulate data}
 \label{table:Exact parameters}
\end{table}

The stress $\sigma (t)$ can be computed exactly using \eqref{eq: sigmaTotal} and is illustrated in figure \ref{abb: strain and stress}.
The goal is to reconstruct the material parameters from table \ref{table:Exact parameters} from given data $\sigma (t)$. All tables present the stiffnesses in megapascals [MPa] and the relaxation times in seconds [s].
 
First, we consider the reconstruction corresponding to the displacement rate of $\dot{\varepsilon}_u=10$ mm/s. The possible number of Maxwell elements is determined by $I:=\lbrace 1, \dots, 5 \rbrace$. Furthermore we assume for the relaxation times $\tau_j\geq \gamma:=0.01$. 
Table \ref{table:non_disturbed_data} shows the results of the cluster algorithm.
While the minimization algorithm determines stiffnesses $\mu_j$ and relaxation times $\tau_j$ for the given maximum number of Maxwell elements (here $N=5$), the cluster algorithm is able to combine them accordingly and reults in the correct number $n^*=3$.

\begin{table}[h]
\centering
		\begin{tabular}{ccccccc}\toprule
			$j=$&$0$ &$1$&$2$&$3$ & $4$ & $5$ \\ \midrule
			Reconstructed values \\ \midrule
			 $\mu_j$ & 10.000 & 4.000 & 3.685 & 1.621 & 1.694 & 1.000 \\ \midrule 
		      $\tau_j$ &-& 0.200 & 3.695 & 3.706 & 3.706 & 25.000 \\\bottomrule
 \toprule
		       After clustering \\ \midrule
			 $\mu_j$ &10.000 & 4.000 & 7.000 & 1.000 \\ \midrule
		 $\tau_j$ & - & 0.200 & 3.700 & 25.000 \\ \bottomrule
		\end{tabular}
  \caption{Reconstructed material parameters before and after clustering}
  \label{table:non_disturbed_data}
\end{table}

While there is only one Maxwell element with $\tau \in [0,1)$ and $\tau \in [10,100)$ after the minimization process, we have three elements in the decade $[1,10)$. The cluster algorithm then combines these elements. Please note, that the parameters of  table \ref{table:Exact parameters} are reconstructed exactly.
The relaxation time $\tau_2=3.695$ is closest to the correct value of $3.7$ before clustering and $\mu_2 \gg \mu_4 > \mu_3$ holds. The weighting \eqref{eq:cluster calculation} then results in these excellent results.

We compare the results of the cluster algorithm with the developed Bayes approach. As mentioned, an advantage of the Bayes algorithm is that it does not require a priori information about the different decades of relaxation times. We choose $q=0.1$ and $M:=\max I=5$. 
Then, similar to the cluster algorithm, we can reconstruct the parameters exactly. The results are listed in table \ref{table:non_disturbed_data Bayes}.

\begin{table}[h]
\centering
		\begin{tabular}{ccccccc}\toprule
		     $j$&$0$ &$1$&$2$&$3$ \\ \midrule
			 $\mu_j$ &10.000 & 4.000 & 7.000 & 1.000 \\ \midrule
		 $\tau_j$ & - & 0.200 & 3.700 & 25.000 \\ \bottomrule
		\end{tabular}
  \caption{Reconstructed material parameters by algorithm \ref{alg:Bayes}}
  \label{table:non_disturbed_data Bayes}
\end{table}

\begin{table}[H]
\centering
		\begin{tabular}{ccccccc}\toprule
		    $j$&$0$ &$1$&$2$&$3$  \\ \midrule
			 $\mu_j$ &9.997 & 57.376 & 7.028 & 0.995 \\ \midrule
		 $\tau_j$ & - & 0.013 & 3.703 & 25.356 \\ \bottomrule
		\end{tabular}
  \caption{Reconstructed parameters by the cluster algorithm for the data from figure \ref{abb: disturbeddata 10}}
  \label{table:disturbeddata 10 cluster}
\end{table}

\begin{table}[H]
\centering
		\begin{tabular}{ccccccc}\toprule
		     $j$&$0$ &$1$&$2$&$3$  \\ \midrule
			 $\mu_j$ &9.9997 & 4.424 & 7.046 & 0.998 \\ \midrule
		 $\tau_j$ & - & 0.162 & 3.693 & 25.307 \\ \bottomrule
		\end{tabular}
  \caption{Reconstructed parameters by algorithm \ref{alg:Bayes} for data from figure \ref{abb: disturbeddata 10}}
  \label{table: disturbeddata 10 Bayes}
\end{table}

\begin{table}[H]
\centering
		\begin{tabular}{ccccccc}\toprule
		     $j$&$0$ &$1$&$2$&$3$  & 4 & 5 \\ \midrule
			 $\mu_j$ &9.9997 & 4.366 & 0.100 & 0.876 & 7.030 & 0.996 \\ \midrule
		 $\tau_j$ & - & 0.100 & 0.100 & 0.339 & 3.701 & 25.347 \\ \bottomrule
		\end{tabular}
  \caption{Reconstructed parameters by minimizing the residual for the data from figure \ref{abb: disturbeddata 10}}
  \label{table:disturbeddata 10 residual}
\end{table}

%%%%%%%%%%%%%%%%%%%%%%

\subsection{Reconstructions from noisy data} \label{sec: disturbed data}

%%%%%%%%%%%%%%%%%%%%%%

We continue by performing reconstructions from noise-contaminated data. We add normally distributed noise to the discretized stress $\sigma$ such that $\lVert \sigma- \sigma^\delta \rVert < \delta$ holds with a noise level $\delta>0$. 
Furthermore, by
\begin{align*}
   \delta_{rel} := \frac{\lVert \sigma- \sigma^\delta \rVert}{\lVert \sigma^\delta \rVert}
\end{align*} 
we denote the relative data error.
In the following experiments the noise level is such that $\delta_{rel} \approx 1 \%$ holds true. 

\begin{figure}[h]
	\centering
		\def\svgwidth{.7\textwidth}
		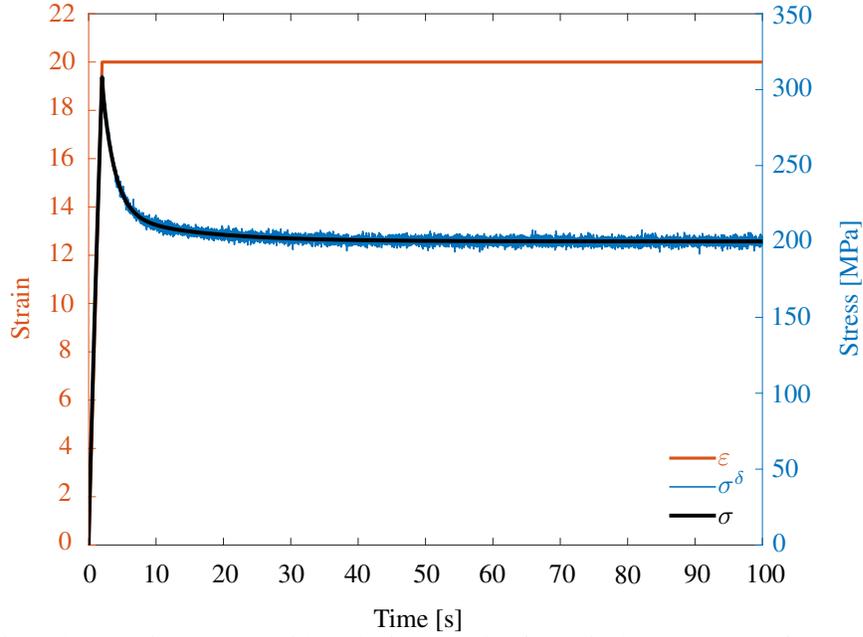
	\caption{Strain and stress-time curves with and without noise for a displacement rate of $\dot{\varepsilon}_u=10$mm/s}
	\label{abb: disturbeddata 10}
\end{figure} 

In figure \ref{abb: disturbeddata 10} both, strain and stress, are plotted. The red line represents the strain curve with a displacement rate of $\dot{\varepsilon}_u=10$ mm/s and a maximum strain of $\bar{\varepsilon}=20\,\%$. The black line is the associated stress with the material parameters from table \ref{table:Exact parameters}. For reconstruction we will use only the noisy stress which is represented by the blue line. 

In table \ref{table:disturbeddata 10 cluster} we see the output of the clustering algorithm. 
Here, first the functional
\begin{align*}
 T_{N}(x) := \lVert \sigma^\delta - F_N (x) \rVert^2
\end{align*} 
is minimized with $N=5$. We realize, that the reconstruction results differ significantly from the exact values. However, it is obvious that the reconstructions of $(\mu_1,\tau_1)$ are the one that are strongest affected by noise, while the other parameters are calculated quite stably. This confirms the considerations of the authors in \cite{Diebels2018}, where it was shown that the reconstruction of small relaxation times is extremely ill-conditioned and this behavior worsens as $\tau \to 0$. Since stiffness and relaxation time are always to be considered in pairs, the deviations in $\tau_1$ also affects the corresponding stiffness $\mu_1$. Since each Maxwell element provides a certain contribution to the total stress, too small values of one parameter are compensated by higher values of the other parameters to match the total stress.

Again, we compare the reconstructions with those from the Bayes method. The results are presented in table \ref{table: disturbeddata 10 Bayes}. Here, in \eqref{eq:Tikhonov_Bayes} we choose $q=0.1$ and $\alpha=1$.
We can see that the results are significantly better than the outcome of the cluster algorithm. The basic stiffness and parameters of the second and third Maxwell elements are reconstructed equally well. However, here we have good approximations of the exact parameters $(4,0.2)$ by the values $(4.4238,0.1617)$ for the first Maxwell element as well. 

The question arises whether this result can also be achieved by simply minimizing the residual with different numbers of $n$. To check this, we implement a third reconstruction method. We set $I=\{n_1,\dots,n_m\}$ and minimize subsequently the residuals
\begin{align} 
\label{eq: residual minimization}
    \lVert \sigma^\delta - F(n_j,x) \rVert^2,\qquad j=1,\ldots,m,
\end{align}
to determine the parameters $x$. Finally we choose $(n,x)$ with minimal residual. 
Since the penalty term, given by the prior, is omitted, the entire statistical inversion aspect is neglected.  

As we can see in table \ref{table:disturbeddata 10 residual}, this reconstruction is very inaccurate. We obtain the maximum possible number of Maxwell elements, however with $\tau_1=\tau_2$. This shows that without additional a priori information on different decades for $\tau_j$ or statistical prior, a reliable reconstruction seems impossible.

\begin{table}[H]
\centering
		\begin{tabular}{lcccc}\toprule
			$j$&$0$ &$1$&$2$ \\ \midrule
			 $\mu^*_j$ [MPa] &5&8&0.5 \\ \midrule 
			 $\tau^*_j$ [s] &-&0.8&50 \\ \bottomrule
		\end{tabular}
	\caption{Second set of material parameters $x^*$ to simulate data}
 \label{table:Exact parameters trial2} 
\end{table}

We want to perform another experiment applying the three methods. To this end, we use the same setting with a maximum strain of $\bar{\varepsilon}=20\,\%$, a displacement rate of $\dot{\varepsilon}_u=10$mm/s and relative noise level of $\delta_{rel} \approx 1 \%$. However, we change the material parameters as listed in table \ref{table:Exact parameters trial2}. Here, the material under consideration has only two Maxwell elements. 

The results of the different algorithms are shown in table \ref{table:disturbeddata 10 trial2 all}. As in the last experiment, the basic stiffness as well as the parameters of the Maxwell element with larger relaxation time are reconstructed by all algorithms in a reasonable accuracy. However, the cluster algorithm as well as the minimization of the residual results in too many Maxwell elements. The cluster algorithm leads to another Maxwell element with very large relaxation time, but very low stiffness that has little influence on the total stress. However, again the first Maxwell element $(\mu^*_1,\tau^*_1)=(8,0.5)$ is reconstructed accurately. Minimizing the residual causes four Maxwell elements where there should be only the first one. Thus, a total of five elements are reconstructed.
Only the Bayes algorithm succeeds in producing a suitable reconstruction of the parameters.

\begin{table}[h]
\centering
		\begin{tabular}{ccccccc}\toprule
            \multicolumn{5}{c}{Cluster algorithm}\\
            \midrule
		     $j$&$0$ &$1$&$2$&$3$  \\
            \midrule
			 $\mu_j$ &4.998 & 29.740 & 0.502 & $8.053\cdot 10^{-10}$ \\ 
            \midrule
		 $\tau_j$ & - & 0.215 & 50.309 & 576.988 \\ 
      \bottomrule
		\end{tabular} \\
\vspace{0.5cm}
\begin{tabular}{ccccccc}\toprule
            \multicolumn{4}{c}{Bayes algorithm}\\
            \midrule
		     $j$&$0$ &$1$&$2$  \\
            \midrule
			 $\mu_j$ &4.998 & 7.923 & 0.502 \\ 
            \midrule
		 $\tau_j$ & - & 0.502 & 50.241  \\ 
      \bottomrule
		\end{tabular} \\
  \vspace{0.5cm}
\begin{tabular}{ccccccc}\toprule
            \multicolumn{7}{c}{Residual minimization \eqref{eq: residual minimization}}\\
            \midrule
		     $j$&$0$ &$1$&$2$&$3$  & 4 & 5  \\
            \midrule
			 $\mu_j$ &4.998 & 0.130 & 0.130 & 4.161 & 3.698 & 0.502 \\ 
            \midrule
		 $\tau_j$ & - & 0.100 & 0.100 & 0.811 & 0.812 & 50.270 \\ 
      \bottomrule
		\end{tabular} \\
  \caption{Reconstructed parameters by the three algorithms to the exact parameters from table \ref{table:Exact parameters trial2}}
  \label{table:disturbeddata 10 trial2 all}
\end{table}

%%%%%%%%%%%%%%%%%%%%%%%%%%%%%%%%%%%%%%%%%%%%%%%%%%

\subsection{Analysis of different displacement rates $\dot{\varepsilon}_u$} \label{sec:different strain rates}

%%%%%%%%%%%%%%%%%%%%%%%%%%%%%%%%%%%%%%%%%%%%%%%%%%%

In the following we consider different displacement rates and how they affect the outcome of the reconstructions. For this purpose, we return to our first experiment corresponding to the parameters from table \ref{table:Exact parameters}. The experimental setup remains the same, the sample is stretched to a maximum of $\bar{\varepsilon}=20\,\%$, the additive noise has a relative noise level of $\delta_{rel}=1\%$. We already know the result to this experiment for a displacement rate of $\dot{\varepsilon}_u=10$ mm/s from table \ref{table: disturbeddata 10 Bayes}. Table \ref{table:disturbeddata 1 Bayes} shows the reconstructions for the Bayes algorithm for this experiment with a displacement rate of $\dot{\varepsilon}_u=1$ mm/s, probability $q=0.1$ and regularization parameter $\alpha=1$.

\begin{table}[h]
\centering
		\begin{tabular}{ccccccc}\toprule
		     $j$&$0$ &$1$&$2$  \\ \midrule
			 $\mu_j$ &9.997 & 7.450 & 1.007  \\ \midrule
		 $\tau_j$ & - & 3.576 & 25.100  \\ \bottomrule
		\end{tabular}
  \caption{Reconstructed parameters using algorithm \ref{alg:Bayes} with $q=0.1$ for $\dot{\varepsilon}_u=1$ mm/s}
  \label{table:disturbeddata 1 Bayes}
\end{table}

As we can see, the Bayes algorithm cannot reconstruct all Maxwell elements reliably, because the first Maxwell element is missing. This is why we want to check, whether the Bayes algorithm is able to find this Maxwell element, if we weaken the strong weighting to a small number of Maxwell elements. For this purpose, table \ref{table:disturbeddata 1 Bayes q03} shows the reconstruction of the Bayes algorithm with $q=0.3$. 

\begin{table}[h]
\centering
		\begin{tabular}{ccccccc}\toprule
		     $j$&$0$ &$1$&$2$& 3  \\ \midrule
			 $\mu_j$ &9.996 & 7.411 & 0.299 & 0.796  \\ \midrule
		 $\tau_j$ & - & 3.500 & 15.369 & 27.100  \\ \bottomrule
		\end{tabular}
  \caption{Reconstructed parameters by algorithm \ref{alg:Bayes} with $q=0.3$ for $\dot{\varepsilon}_u=1$ mm/s}
  \label{table:disturbeddata 1 Bayes q03}
\end{table}

In fact, three Maxwell elements are reconstructed this time, but the first Maxwell element $(\mu^*_1,\tau^*_1)=(4,0.2)$ is not determined. Moreover, the third Maxwell element $(\mu^*_3,\tau^*_3)=(1,25)$ is split into two Maxwell elements $(0.2988,15.3686)$ and $(0.7955,27.1003)$.\\

\begin{figure}
     \centering
    \subfloat{\def\svgwidth{0.5\textwidth}
     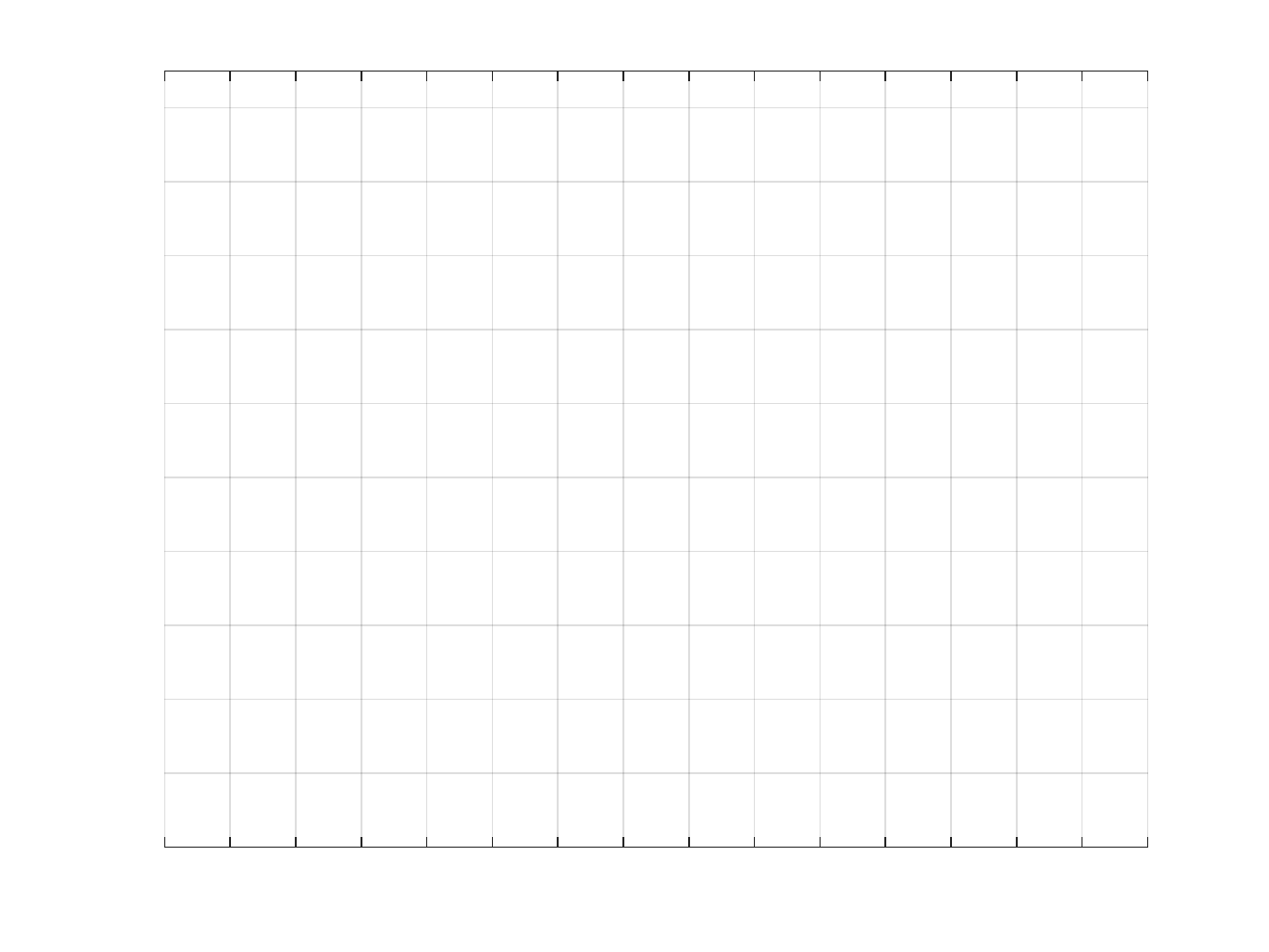}
 \subfloat{\def\svgwidth{0.5\textwidth}
   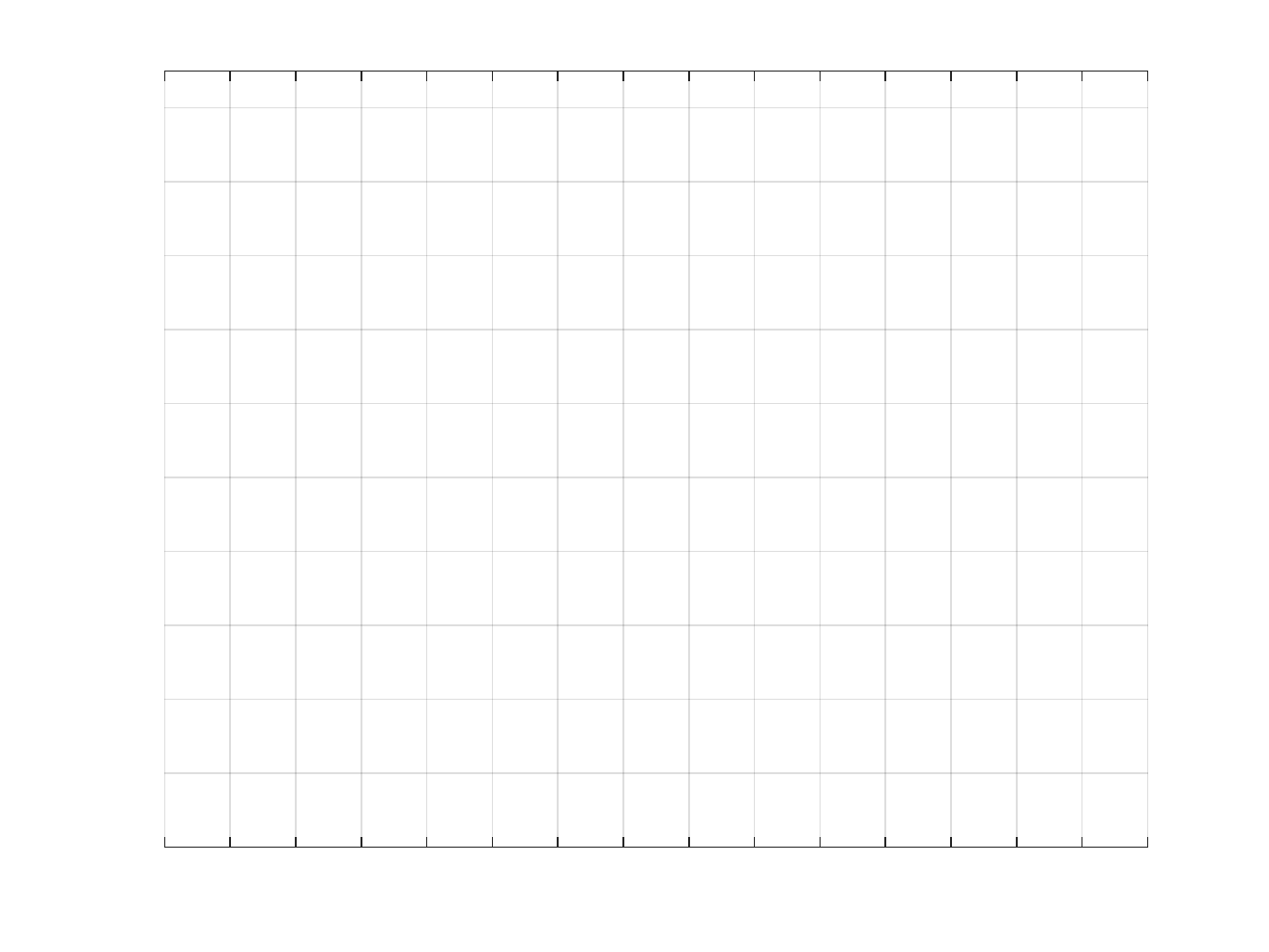}
   \caption{Individual stress components for a displacement rate of $1$ mm/s and $10$ mm/s}
        \label{abb:displacement rates}
\end{figure}

The displacement rate seems somehow to affect the reconstruction outcome. To this end we consider the different stress-time curves generated by a displacement rate of $1$ mm/s and $10$ mm/s (Figure \ref{abb:displacement rates}). For each of these curves, the contribution of each element is demonstrated separately. 
That is, $\sigma_0$ is the stress generated by the single spring, while $\sigma_{1,2,3}$ is the stress of the three Maxwell elements. The sum of these individual stresses yields the total stress \eqref{eq: sigmaTotal}.
The maximum strain is $20\,\%$, which is achieved at a displacement rate of $1$mm/s and $10$ mm/s after $20$ seconds and $2$ seconds, respectively. This represents a significant difference for the individual stresses of the Maxwell elements. In both cases, the maximum stress of the single spring is $200$ MPa, but at different time instances. 
However, the maximum stress values of the Maxwell elements are much lower at the slower displacement rate of $1$ mm/s than at a higher rate. In both cases, the maximum is achieved at $t=\bar{\varepsilon}/\dot{\varepsilon}$. Table \ref{table: maxStress ME different strain rates} lists the different values. 
For a slower displacement rate, the stress values of the Maxwell elements are non-zero over a longer period of time, recognizable by the second Maxwell element in figure \ref{abb:displacement rates}. However, since the values especially for the first Maxwell element are much smaller than at a higher displacement rate, these small values are much more affected by noise. 
Therefore, it is advisable to use higher shift rates. For this reason, we confine ourselves to a displacement rate of $10$mm/s in the next examples.

\begin{table}	[h]
	\centering
	\begin{tabular}{c|c|ccc} \toprule &	&$j=1$	&$j=2$	&$j=3$ 			\\ \hline
	$\dot{\varepsilon}_u = 1$ 	&	$\sigma_j(\bar{\varepsilon}/\dot{\varepsilon})$				&0.4	&12.94	&9.9 		\\
				 \hline
		
	$\dot{\varepsilon}_u = 10$	&	$\sigma_j(\bar{\varepsilon}/\dot{\varepsilon})$				&4		&85.57	&18.48	 	\\ \bottomrule
	\end{tabular}
 \caption{Maximum values of individual stress components with slow and fast displacement rates at time $t=\bar{\varepsilon}/\dot{\varepsilon}$}
		\label{table: maxStress ME different strain rates}
\end{table}

%%%%%%%%%%%%%%%%%%%%%%%%%%%%%%%%%%%%%%%%%%%%%%%%%%

\subsection{Effect of $q$ on the reconstruction results}
\label{sec:num_effect_q}

%%%%%%%%%%%%%%%%%%%%%%%%%%%%%%%%%%%%%%%%%%%%%%%%%&

As already deduced from table \ref{table:disturbeddata 1 Bayes q03}, increasing $q$ causes the algorithm to favor a higher number of Maxwell elements as intended by the prior (c.f. section \ref{sec: effect q prior}). In this subsection we investigate the sensitivity of the results regarding $q$. 

We consider again the second experiment as already described in section \ref{sec: disturbed data} with material parameters as listed in table \ref{table:Exact parameters trial2}. We already know that the Bayes algorithm for $q=0.1$ reconstructs the parameters very well. We change now the value of $q$ to $q=0.9$. The prior $\rho_0(n)$, thus, favors a higher number of Maxwell elements, that is, the penalty term $\phi_0(n)$ is monotonically decreasing and attains its maximum for $n=1$.

\begin{table}[h]
\centering
		\begin{tabular}{ccccccc}\toprule
		     $j$&$0$ &$1$&$2$& 3  & 4 & 5 \\ \midrule
			 $\mu_j$ & 4.998 & 0.130 & 0.130 & 4.161 & 3.698 & 0.502\\ \midrule
		 $\tau_j$ & - & 0.100 & 0.100 & 0.811 & 0.812 & 50.270 \\ \bottomrule
		\end{tabular}
 \caption{Reconstructed parameters of the second experiment by algorithm \ref{alg:Bayes} with $q=0.9$}
\label{table: disturbeddata2 Bayes q09}
\end{table}

Table \ref{table: disturbeddata2 Bayes q09} shows the reconstructed parameters for this experiment. As we can see, the Bayes algorithm reconstructs five elements instead of two. This is consistent with the prior. %Additionally, since it is easier to minimize the residual $\lVert \sigma^\delta - F(n,x) \rVert$ for a higher number of parameters, a small number of Maxwell elements is implausible. 
A test series using different values of success probability $q$ increasing with step size $0.1$ reveals, that only for $q=0.4$ the correct number of Maxwell elements is reconstructed, see table \ref{table: disturbeddata2 Bayes q04}. Further tests prove that for $0< q \leq 0.4$ the material parameters are reliably reconstructed. This includes the parameters that we obtained for $q=0.1$ (c.f. \ref{table:disturbeddata 10 trial2 all}). 

\begin{table}[h]
\centering
		\begin{tabular}{ccccccc}\toprule
		     $j$&$0$ &$1$&$2$\\ \midrule
			 $\mu_j$ & 4.998 & 7.923 & 0.502  \\ \midrule
		 $\tau_j$ & - & 0.807 & 50.241 \\ \bottomrule
		\end{tabular}
 \caption{Reconstructed parameters of the second experiment by algorithm \ref{alg:Bayes} with $q=0.4$}
\label{table: disturbeddata2 Bayes q04}
\end{table}

A third trial uses the same setting with $\bar{\varepsilon}=20\%$, $\dot{\varepsilon}_u=10$ mm/s, and $\delta_{rel} \approx 1\%$. The parameters are modified according to table \ref{table:exact parameters trial 3}. The experiment is extended to $T=10 000$ seconds, since the largest relaxation time is $\tau_5=1200$ s. It is worth mentioning that in the previous experiments we used a temporal sampling rate of $\Delta t=0.01$ seconds. That is, $t_i=i \cdot \Delta t$ for $i=0,\ldots,m$. In the third experiment we choose the same $\Delta t=0.1$ s to obtain a sampling rate which is small enough to catch the influence of the first Maxwell element having a relaxation time of $\tau_1=0.8$ s. Moreover we choose a probability of $q=0.1$ which favors a small number of elements. Table \ref{table:perturbeddata3 Bayes q01} presents the reconstruction result.

\begin{table}[h]
\centering
		\begin{tabular}{ccccccc}
            \toprule
		     $j$&$0$ &1 & 2 & 3 & 4 & 5    \\ 
       \midrule
			 $\mu_j^*$ & 10 & 8 & 7 & 1 & 4 & 0.5 \\ 
            \midrule
		 $\tau_j^*$ & - & 0.8 & 3.7 & 25 & 500 & 1200 \\ 
            \bottomrule
		\end{tabular}
 \caption{Exact parameters $x^*$ for the third experiment}
  \label{table:exact parameters trial 3}
\end{table}

\begin{table}[h]
\centering
		\begin{tabular}{ccccccc}
  \toprule
		     $j$&$0$ &$1$&$2$ & 3 & 4 & 5     \\ 
       \midrule
			 $\mu_j$ & 10.000 & 8.427 & 7.011 & 0.893 & 4.077 & 0.419 \\ 
    \midrule
		 $\tau_j$ & - & 0.763 & 3.846 & 27.47 & 505.8 & 1282.5 \\ 
   \bottomrule
		\end{tabular}
 \caption{Reconstructed parameters of the third experiment by the Bayes algorithm with $q=0.1$}
  \label{table:perturbeddata3 Bayes q01}
\end{table}

Although the algorithm prefers small numbers of Maxwell elements, all elements and the corresponding parameters are reconstructed sufficiently accurate. This results suggests to use a small $q$, like, e.g., $q=0.1$, if the number of Maxwell elements is unknown, since the algorithm is able to deliver the correct number as output.
%This again suggests that it is not possible to describe the data term with a lower number of Maxwell elements and thus parameters. In any case, the noisy stress curve is more easily approximated by a higher number of parameters. Therefore, it is advisable to choose a low success probability $q$, so that the algorithm prefers a low number of Maxwell elements. Material models with a high number of Maxwell elements can also be reconstructed with it. 

%%%%%%%%%%%%%%%%%%%%%%%%%%%%%%%%%%%%%%%%%%%%%%%%%%

\subsection{Regularization by an additional penalty term} \label{sec: regularization}

%%%%%%%%%%%%%%%%%%%%%%%%%%%%%%%%%%%%%%%%%%%%%%%%%%

As outlined in section \ref{chap: Bayes Inversion}, penalty terms are necessary for a stable reconstruction of the parameters. Regularization was also applied to the clustering algorithm in the presence of noise in the data, c.f. \cite{Rothermel2022}. The numerical results emphasize that the penalty term leads to more accurate results as the cluster algorithm without penalty.
The only data sets, which are not recovered accurately, are those with a displacement rate of $\dot{\varepsilon}_u=1$ mm/s. Although, this can be avoided by a corresponding experimental setup, we investigate, whether an additional penalty term leads to better reconstruction results. 
As penalty term we use 
\begin{align} 
    \Omega_1(x)=\frac{1}{2}\lVert x\rVert ^2.
    \label{eq: penaltyTikh}
\end{align}
For fixed number of Maxwell elements $n$, this leads to a classical Tikhonov-Phillips regularization of the form
\begin{align}
       \min_{x \in \mathcal{D}(F_n)} T_{\alpha,n}(x):= \min_{x \in \mathcal{D}(F_n)} \left\{ \frac{1}{2} \lVert F_n(x) - \sigma^\delta\rVert^2 +\alpha \Omega_1(x) \right\},
\end{align} 
see, e.g., \cite{Engl1989ConvergenceProblems, Kaltenbacher2008IterativeProblems, Neubauer1992TikhonovScales,Scherzer1993TheProblems}. The parameter $\alpha>0$ acts as a regularization parameter and balances the influence of the data term and the penalty term on the minimizer.

\begin{table}[h]
\centering
		\begin{tabular}{ccccccc}
  \toprule
		     $j$&$0$ &$1$&$2$ & 3 &    \\ 
       \midrule
			 $\mu_j$ & 10.000 & 3.765 & 3.765 & 1.117 \\ 
    \midrule
		 $\tau_j$ & - & 3.376 & 3.376 & 23.238 \\ 
   \bottomrule
		\end{tabular}
 \caption{Reconstructed parameters of the first experiment for a displacement rate $\dot{\varepsilon}_u=1$ mm/s by the Bayes algorithm with additional penalty term \eqref{eq: penaltyTikh} and $\alpha=0.5$}
  \label{table:disturbeddata 1mm/s Bayes Tikh}
\end{table}

Table \ref{table:disturbeddata 1mm/s Bayes Tikh} shows the result of the Bayes algorithm with penalty \eqref{eq: penaltyTikh} and regularization parameter $\alpha=0.5$. We see that the algorithm now reconstructs three Maxwell elements, however the values of the first and second Maxwell element are identical and the correct first Maxwell element $(\mu_1^*, \tau_1^*)=(4,0.2)$ is not recovered at all. As expected, the penalty term here ensures that large values are penalized. This also explains the reduction of the stiffness value $\mu_3$. As a remedy, we apply a further penalty term of the form 
\begin{align} \label{eq: penalty term tau}
    \Omega_2(x)= \frac{1}{2} \tau_1^2
\end{align}
to penalize large values $\tau_1$. The regularization parameter is increased to $\alpha=100$. 

\begin{table}[h]
\centering
		\begin{tabular}{ccccccc}
  \toprule
		     $j$&$0$ &$1$&$2$ & 3 &    \\ 
       \midrule
			 $\mu_j$ & 9.997 & 2.015 & 7.327 & 1.001 \\ 
    \midrule
		 $\tau_j$ & - & 0.100 & 3.621 & 25.190 \\ 
   \bottomrule
		\end{tabular}
  \caption{Reconstructed parameters of the first experiment for a displacement rate $\dot{\varepsilon}_u=1$ mm/s by the Bayes algorithm with additional penalty term \eqref{eq: penalty term tau} and $\alpha=100$}
  \label{table:disturbeddata 1mm/s Bayes tau}
\end{table}

Table \ref{table:disturbeddata 1mm/s Bayes tau} shows the corresponding outcome of this experiment. We see a significant improvement compared to the results of table \ref{table:disturbeddata 1 Bayes} and table \ref{table:disturbeddata 1mm/s Bayes Tikh}. All three Maxwell elements are reconstructed accurately and the third Maxwell element is no longer negatively affected by $\Omega_2$, though it does only take $\tau_1$ into account. We admit that the first Maxwell element still shows some error sensitivity, even though the reconstruction means a significant improvement.

We conclude that the recommendation to use a higher displacement rate, such as $\dot{\varepsilon}_u=10$ mm/s, remains valid even with an additional penalty term.

%%%%%%%%%%%%%%%%%%%%%%%%%%%%%%%%%%%%%
\section{Conclusion}
%%%%%%%%%%%%%%%%%%%%%%%%%%%%%%%%%%%%%

In this article, we considered the inverse problem of identifying material parameters in a viscoelastic structure using a generalized Maxwell model. One major challenge is the fact that the number of Maxwell elements in this model, and, thus, the number of material parameters, are unknown. Based on statistical inversion using a binomial prior we developed a novel reconstruction method which is able to compute the number of elements along with the corresponding material parameters in a stable way. Since the forward operator acts on $\mathbb{N}\times \l2$, it was as a further novelty necessary to extend the existing regularization theory to semigroups, where we equipped $\mathbb{N}$ with the discrete topology. The influence of the success probability $q$ has been studied and the method is compared to the cluster algorithm by extensive numerical tests. While the cluster algorithm is very sensitive with respect to noise in data, the Bayes algorithm proves to be stable, where in case of low displacement rates additional penalty terms improve the reconstruction results. Using statistical inversion theory yields a penalty term that is tailored to the problem and does not demand for any a priori information on relaxation times and number of Maxwell elements. Controling the success probability $q$ allows to adapt the reconstruction according to any prior information on different structures. In general, it is advisable to choose a small value $q$ in order to favor a small number of Maxwell elements, whereas also larger $n$ can be reconstructed reliably. The numerical evaluation proved the stability of the new method for different paramter settings as well as its superiority with respect to the cluster algorithm.

%%%%%%%%%%%%%%%%%%%%%%%%%%%%%%%%%%%%

\bibliography{references} 

%%%%%%%%%%%%%%%%%%%%%%%%%%%%%%%%%%%%

\end{document}